\documentclass[11pt]{article}
\usepackage{amsfonts,latexsym,amsmath,amscd,geometry}
\usepackage{amssymb}
\usepackage{latexsym}

\newcommand \nc{\newcommand}
\newtheorem{theorem}{Theorem}[section]
\newtheorem{lemma}[theorem]{Lemma}

\newtheorem{definition}[theorem]{Definition}

\nc{\ba}{\begin{array}}\nc{\ea}{\end{array}}
\nc{\be}{\begin{eqnarray}}\nc{\ee}{\end{eqnarray}}
\nc{\beq}{\begin{equation}}\nc{\eeq}{\end{equation}}
\nc{\bex}{\begin{eqnarray*}}\nc{\eex}{\end{eqnarray*}}
\nc{\btm}{\begin{theorem}} \nc{\etm}{\end{theorem}}
\nc{\blm}{\begin{lemma}} \nc{\elm}{\end{lemma}}
\nc{\R}{\mathbb{R}} \nc{\va}{\varepsilon} \nc{\ls}{\limits}

\def\pf{\noindent{\bf Proof.\quad}}

\newcommand \qed {\hfill $\Box$}

\begin{document}

\thispagestyle{empty}
\title{Energy identity of approximate biharmonic maps to Riemannian manifolds and its application}
\author{
Changyou Wang\footnote{ Department of Mathematics, University of Kentucky, Lexington, KY 40506,
cywang@ms.uky.edu.}\ \quad
Shenzhou Zheng\footnote{Department of Mathematics, Beijing Jiaotong University, Beijing 100044, P. R. China, shzhzheng@bjtu.edu.cn.}}
\maketitle

\begin{abstract} We consider in dimension four
weakly convergent sequences of approximate biharmonic maps to a Riemannian manifold with bi-tension fields bounded in $L^p$ for $p>\frac43$.
We prove an energy identity that accounts for the loss of hessian energies by the sum of hessian energies over finitely many nontrivial biharmonic maps on $\mathbb R^4$.
As a corollary, we obtain an energy identity for the heat flow of biharmonic maps at time
infinity.
\end{abstract}

\section {Introduction}

This is a continuation of our previous work \cite{WZ} on the blow-up analysis of approximate biharmonic maps in dimension 4. In \cite{WZ}, we obtained in dimension four 
an energy identity for approximate biharmonic maps into sphere with bounded $L^p$ bi-tension field for $p>1$, and an energy identity of the heat flow of biharmonic map into sphere at time infinity.
The aim of this paper is to extend the main theorems of \cite{WZ} to any compact Riemannian manifold
without boundary, under the additional assumption that the bi-tension
fields are bounded in $L^p$ for $p>\frac43$. The main results of this paper was announced in \cite{WZ}. 

Let $\Omega\subset \mathbb{R}^4$ be a bounded smooth domain, and $(N,h)$ be
a compact $n$-dimensional Riemannian manifold without boundary, embedded into $K$-dimensional Euclidean space $\mathbb R^{K}$.
Recall the Sobolev space $W^{l,p}(\Omega,N)$,
$1\le l<+\infty$ and $1\le p<+\infty$, is defined by
$$
W^{l,p}(\Omega, N)
=\Big\{v\in W^{l,p}(\Omega, \mathbb R^{K}):
\ v(x)\in N \ {\rm{\ a.e.\ }}x\in\Omega\Big\}.
$$

In this paper we will discuss the limiting behavior of weakly convergent sequences of approximate
(extrinsic) biharmonic maps $\displaystyle \{u_k\}\subset W^{2,2}(\Omega,N)$ in dimension $n=4$, especially an energy identity
during the process of convergence.
First we recall the notion of approximate (extrinsic) biharmonic maps.
\begin{definition} A map $u\in W^{2,2}(\Omega,N)$ is called an approximate biharmonic map if there exists a bi-tension
field $h\in L^1_{\rm{loc}}(\Omega, \mathbb R^{K})$ such that
\begin{equation}\label{approx_biharm}
\Delta^2 u=\Delta(\mathbb B(u)(\nabla u,\nabla u))+2\nabla\cdot\langle\Delta u,
\nabla(\mathbb P(u))\rangle-\langle\Delta(\mathbb P(u)),\Delta u\rangle+h
\end{equation}
in the distribution sense, where $\mathbb P(y):\mathbb R^{K}\to T_y N$ is
the orthogonal projection from $\mathbb R^{K}$ to the tangent space of $N$ at $y\in N$,
and $\displaystyle\mathbb B(y)(X,Y)=-\nabla_X\mathbb P(y)(Y), \ \forall  X, Y\in T_y N$, is the second fundamental form of $N\subset\mathbb R^{K}$.
In particular, if $h=0$ then $u$ is called a biharmonic map to  $N$.
\end{definition}

Note that biharmonic maps to a Riemannian manifold $N$ are critical points of the hessian energy functional
$$\displaystyle E_2(u)=\int_\Omega |\nabla^2 u|^2\,dx$$
over $W^{2,2}(\Omega, N)$.
Biharmonic maps are  higher order extensions of harmonic maps.
The study of regularity of biharmonic maps has generated considerable interests after the initial work by Chang-Wang-Yang \cite{CWY},  the readers
can refer to Wang \cite{W1,W2,W3}, Strzelecki \cite{S},  Lamm-Rivi\`ere \cite{LR},
Struwe \cite{Struwe},  Scheven \cite{scheven1, scheven2} (see also Ku \cite{K} and Gong-Lamm-Wang \cite{GLW} for the boundary regularity).  In particular,  the interior regularity theorem asserts that the smoothness of $W^{2,2}$-biharmonic maps holds in dimension $n=4$,
and the partial regularity of stationary $W^{2,2}$-biharmonic maps holds in dimensions $n\ge 5$.

It is an important observation that biharmonic maps are invariant under dilations in $\mathbb R^n$ for $n=4$. Such
a property leads to non-compactness of biharmonic maps in dimension $4$, which prompts recent studies
by Wang \cite{W1} and Hornung-Moser \cite{HM} concerning the failure of strong convergence for weakly convergent biharmonic maps.  Roughly speaking,  the results in \cite{W1} and \cite{HM} assert that the failure of strong convergence occurs at finitely many concentration points of hessian energy, where finitely many
bubbles (i.e. nontrivial biharmonic maps on $\mathbb R^4$) are generated, and the total hessian energies
from these bubbles account for the total loss of hessian energies during the process of convergence.

Our first result is to extend the results from \cite{W1} and \cite{HM} to the context of suitable approximate biharmonic maps to a compact Riemannian manifold $N$. More precisely, we have
\begin{theorem}\label{energy_identity} For $n=4$, suppose $\{u_k\}\subset W^{2,2}(\Omega,N)$ is a sequence of approximate biharmonic maps, which are bounded in $W^{2,2}(\Omega, N)$ and have their bi-tension fields
$h_k$ bounded in $L^p$ for $p>\frac 43 $, i.e.
\begin{equation}\label{bound}
M:=\sup_{k}\Big(\|u_k\|_{W^{2,2}}+\|h_k\|_{L^p}\Big)<+\infty.
\end{equation}
Assume $u_k\rightharpoonup u$ in $W^{2,2}$ and $h_k\rightharpoonup h$ in $L^p$. Then \\
(i) $u$ is an approximate biharmonic map to  $N$ with $h$ as its bi-tension field. \\
(ii) there exist a nonnegative integer $L$ depending on $M$
and $L$ points $\{x_1,\cdots,x_L\}\subset\Omega$
such that
$\displaystyle u_k\rightarrow u \ strongly\ in\ 
W^{2,2}_{\rm{loc}}\cap C^0_{\rm{loc}}(\Omega\backslash\{x_1,\cdots,x_L\},N).$\\
(iii) For $1\le i\le L$, there exist a positive integer $L_i$ depending on $M$ and
$L_i$ nontrivial smooth biharmonic map
$\omega_{ij}$ from $\mathbb R^4$ to $N$ with finite hessian energy, $1\le j\le L_i$, such that
\begin{equation}\label{energy_id1}
\lim_{k\rightarrow\infty}\int_{B_{r_i}(x_i)} |\nabla^2u_k|^2
=\int_{B_{r_i}(x_i)}|\nabla^2 u|^2+\sum_{j=1}^{L_i}\int_{\mathbb{R}^4}|\nabla^2\omega_{ij}|^2,
\end{equation}
and
\begin{equation}\label{energy_id2}
\lim_{k\rightarrow\infty}\int_{B_{r_i}(x_i)} |\nabla u_k|^4
=\int_{B_{r_i}(x_i)} |\nabla u|^4
+\sum_{j=1}^{L_i}\int_{\mathbb{R}^4} |\nabla\omega_{ij}|^4,
\end{equation}
where $\displaystyle
r_i=\frac 12\min_{1\le j\le L,\ j\neq
i}\left\{|x_i-x_j|,\ {\rm{dist}}(x_i,\partial\Omega)\right \}.$
\end{theorem}




As an application of Theorem 1.2, we study asymptotic behavior at time infinity for the heat flow of biharmonic maps in dimension $4$.

Let's review the  studies on the heat flow of biharmonic maps undertaken by Lamm \cite{Lamm}, Gastel \cite{Gastel}, Wang \cite{W4}, and Moser \cite{M}.  
The equation of heat flow of (extrinsic) biharmonic maps into $N$ is to seek
$u:\Omega\times [0,+\infty) \to N$ that solves:
\begin{eqnarray}\label{heat_biharm}
u_t+\Delta^2 u&=&\Delta(\mathbb B(u)(\nabla u,\nabla u))+2\nabla\cdot\langle\Delta u,
\nabla(\mathbb P(u))\rangle-\langle\Delta(\mathbb P(u)),\Delta u\rangle,
\ \Omega\times (0,+\infty) \\
u&=& u_0, \ \Omega\times \{0\}\label{init_cond}\\
(u,\frac{\partial u}{\partial\nu})&=&(u_0,
\frac{\partial u_0}{\partial\nu}),
\ \partial\Omega\times (0,+\infty),\label{bdry_cond}
\end{eqnarray}
where $u_0\in W^{2,2}(\Omega, N)$ is a given map. Note that any time independent solution
$u:\Omega\to N$ of (\ref{heat_biharm}) is a biharmonic map to $N$.

In dimension $n=4$, Lamm \cite{Lamm} established the existence of global smooth solutions to (\ref{heat_biharm})-(\ref{bdry_cond}) for $u_0\in W^{2,2}(\Omega, N)$ with small $W^{2,2}$-norm, and Gastel \cite{Gastel} and Wang \cite{W4} independently showed that
there exists a unique global weak solution to (\ref{heat_biharm}))-(\ref{bdry_cond}) for any initial data $u_0\in W^{2,2}(\Omega, N)$ that has at most finitely many singular
times. Moreover, such a solution enjoys the energy inequality:
\begin{equation}\label{energy_ineq}
2\int_0^T\int_\Omega |u_t|^2+\int_\Omega |\Delta u|^2(T)
\le \int_\Omega |\Delta u_0|^2,
\ \forall \ 0<T<+\infty.
\end{equation}
Recently, Moser \cite{M} showed the existence of a global weak solution to (\ref{heat_biharm})-(\ref{bdry_cond})
for any target manifold $N$ in dimensions $n\le 8$.

It follows from (\ref{energy_ineq}) that there exists a sequence $t_k\uparrow\infty$ such that $u_k:=u(\cdot, t_k)\in W^{2,2}(\Omega,N)$ satisfies \\
(i) $\displaystyle \tau_2(u_k):=u_t(t_k)$ satisfies $\|\tau_2(u_k)\|_{L^2}\rightarrow 0$; and\\
(ii) $u_k$ satisfies in the distribution sense
\begin{equation}\label{approx_biharm1}
-\Delta^2 u_k+\Delta(\mathbb B(u_k)(\nabla u_k,\nabla u_k))+2\nabla\cdot\langle\Delta u_k,\nabla(\mathbb P(u_k))\rangle-\langle\Delta(\mathbb P(u_k)),\Delta u_k\rangle=\tau_2(u_k).
\end{equation}
By Definition 1.1 $\{u_k\}$ is a  sequence of approximate biharmonic maps to $N$, which are bounded in $W^{2,2}$ and have their bi-tension fields bounded in $L^2$. Hence,
as an immediate corollary, we obtain 
\begin{theorem}\label{energy_identity1} For $n=4$ and $u_0\in W^{2,2}(\Omega, N)$, let
$u:\Omega\times \mathbb R_+\to N$, with $u\in L^\infty(\mathbb R_+, W^{2,2}(\Omega))$ and
$u_t\in L^2(\mathbb R_+, L^2(\Omega))$, be a global weak solution of (\ref{heat_biharm})-(\ref{bdry_cond}) that satisfies the
energy inequality (\ref{energy_ineq}). Then there exist $t_k\uparrow +\infty$, a biharmonic map $u_\infty\in C^\infty\cap W^{2,2}(\Omega,N)$
with $u_\infty=u_0$ on $\partial\Omega$, and a nonnegative integer $L$ and  $L$ points
$\{x_1,\cdots, x_L\}\subset\Omega$ such that\\ (i) $u_k:=u(\cdot,t_k)\rightharpoonup u_\infty$ in $W^{2,2}(\Omega, N)$.\\
(ii) $u_k\rightarrow u_\infty$ in $C^0_{\rm{loc}}\cap W^{2,2}_{\rm{loc}}(\Omega\setminus\{x_1,\cdots, x_L\}, N)$.\\
(iii) for $1\le i\le L$, there exist a positive integer $L_i$ and $L_i$ nontrivial biharmonic maps $\{\omega_{ij}\}_{j=1}^{L_i}$
on $\mathbb R^4$ with finite hessian energies such that
\begin{equation}\label{energy_id1}
\lim_{k\rightarrow\infty}\int_{B_{r_i}(x_i)} |\nabla^2u_k|^2
=\int_{B_{r_i}(x_i)}|\nabla^2 u_\infty|^2+\sum_{j=1}^{L_i}\int_{\mathbb{R}^4}|\nabla^2\omega_{ij}|^2,
\end{equation}
and
\begin{equation}\label{energy_id2}
\lim_{k\rightarrow\infty}\int_{B_{r_i}(x_i)} |\nabla u_k|^4
=\int_{B_{r_i}(x_i)} |\nabla u_\infty|^4
+\sum_{j=1}^{L_i}\int_{\mathbb{R}^4} |\nabla\omega_{ij}|^4,
\end{equation}
where $\displaystyle
r_i=\frac 12\min_{1\le j\le L,\ j\neq
i}\left\{|x_i-x_j|,\ {\rm{dist}}(x_i,\partial\Omega) \right\}.$
\end{theorem}


The main ideas to prove Theorem 1.2 can be outlined as follows. First, we adapt the arguments from \cite {W2, W3} to establish an $\varepsilon_0-$regularity theorem for any approximate biharmonic map $u$ with bi-tension field $h\in L^p$ for
$p>1$, which asserts that $u\in C^\alpha$ for any $\alpha\in (0,\min\{1,4(1-\frac{1}{p})\})$
and $\nabla^4 u\in L^p$. Second, we prove that for $p>\frac43$ there is no concentration of  
angular hessian energy in the neck region by comparing the approximate biharmonic maps with
radial biharmonic functions over annulus. This is a well known technique in harmonic maps
in dimension two (see \cite{DT}). For biharmonic maps in dimension 4, it was derived by
Hornung-Moser \cite{HM}.  Third, we use a Pohozaev type argument to control the radial hessian
energy by the angular hessian energy and $L^p$-norm of bi-tension fields in the neck region.
The assumption $p\ge\frac43$ seems to be necessary to validate the Pohozaev type argument, since we need
$\Delta^2 u_k\cdot(x\cdot\nabla u_k) \in L^1$ and $h_k\cdot(x\cdot\nabla u_k)\in L^1$. 
It remains to be an open question whether Theorem 1.2 holds for $1<p\le\frac43$.

The paper is organized as follows. In \S 2, we establish the H\"older continuity and
$W^{4,p}$-regularity for any approximate biharmonic map with its bi-tension field in $L^p$ for 
$p>1$. In \S3, we show the strong convergence under  the smallness condition
of hessian energy and set up the bubbling process. In \S4, we prove Theorem 1.2 by
establishing (i) there is no concentration of angular hessian energy
in the neck region; and (ii) control the radial hessian energy in the neck region by
angular hessian energy and $L^p$-norm of bi-tension field through a Pohozaev type argument.

\bigskip
\noindent{\bf Acknowledgement}.\  The first author is partially supported by NSF 1000115. The second author is partially supported by NSFC grant 11071012.
This work is conducted while the second author is visiting the University of Kentucky, he would like to thank Department of Mathematics for its support and hospitality.

\bigskip
\noindent{\bf Added Note}.  After we completed this paper, we noticed from a closely related preprint posted  in arxiv.org at December 23, 2011 by Laurain and Rivi\`ere (arXiv:1112.5393v1), in which they claimed Theorem 1.2 holds for all $p>1$.  Since the main ideas of our proof  are different from theirs,
we believe that our work here shall have its own interest.

\section{A priori estimates of approximate biharmonic maps}
\setcounter{equation}{0}
\setcounter{theorem}{0}

In this section, we will establish both  H\"older continuity and 
$W^{4,p}$-regularity  for any approximate biharmonic map with bi-tension field 
$h\in L^p$ for $p>1$, under the smallness condition of hessian energy. The proof of H\"older continuity is based
on suitable modifications and extensions of that by Wang \cite{W2} Theorem A on the regularity of
biharmonic maps to the context of approximate biharmonic maps in dimension 4. The proof of
a bootstrap type argument, which may have its own interest. It is needed for the $\varepsilon_0$-compactness lemma and 
Pohozaev argument for approximate biharmonic maps. 

Denote by $B_r(x)\subset\mathbb R^4$ the ball with center $x$ and radius $r$, and $B_r=B_r(0)$.  First, we have
\begin{lemma} \label{holder} For any $\alpha\in \Big(0, \min\{1, 4(1-\frac{1}{p})\}\Big)$, 
there exist $\varepsilon_0>0$ such that if $u\in W^{2,2}(B_2,N)$ is an approximate biharmonic map with its bi-tension field $h\in L^p(B_1)$ for $p>1$, and satisfies
\begin{equation}\label{smallness1}
\int_{B_2}\Big(|\nabla^2 u|^2+ |\nabla u|^4\Big) \le\varepsilon^2_0.
\end{equation}
then $u\in C^{\alpha}(B_{\frac 12},N)$ and
\begin{equation}\label{holder_est}
\Big[u\Big]_{C^{\alpha}(B_{\frac 12})}
\le C\Big(\varepsilon_0+\|h\|_{L^p(B_1)}\Big).
\end{equation}
\end{lemma}
\pf We follow \cite{W2} \S2, \S3 and \S4 closely and only sketch the main steps of proof here. 
First, by choosing $\varepsilon_0$ sufficiently small, Proposition 3.2 and Theorem 3.3 of \cite{W2}
imply that there exists an adopted frame $\{e_\alpha\}_{\alpha=1}^m$ ($m={\rm{dim}} N$)
along with $u^*TN$ on $B_1$ such that its connection form $A=(\langle De_\alpha, e_\beta\rangle)$ satisfies:
\begin{equation}
\begin{cases}
d^*A=0 \ {\rm{in}}\ B_1; \ \ x\cdot A=0 \ {\rm{on}}\ \partial B_1\\
\|A\|_{L^4(B_1)}+\|\nabla A\|_{L^2(B_1)} \le C\|\nabla u\|_{L^4(B_1)}^2\\
\|\nabla A\|_{L^{2,1}(B_\frac12)}+\|A\|_{L^{4,1}(B_\frac12)}\le C(\|\nabla^2 u\|_{L^2(B_1)}+\|\nabla u\|_{L^4(B_1)}).
\end{cases}
\end{equation}
Here $L^{r,s}$ denotes Lorentz spaces for $1\le r<+\infty, 1\le s\le\infty$, see \cite{W2} \S2 or \cite{H} for its definition and basic properties. 

Utilizing $\{e_\alpha\}_{\alpha=1}^m$, we can rewrite (\ref{approx_biharm}) into the following form (see \cite{W2} Lemma 4.1):
for $1\le \alpha \le m$, 
\begin{equation}\label{approx-express}
\Delta\nabla\cdot\langle\nabla u, e_{\alpha}\rangle=G(u,e_{\alpha})+\langle h, e_{\alpha}\rangle,
\end{equation}
where
\begin{equation}\label{approx-detail}
G(u, e_{\alpha})=\begin{cases}\Delta\langle\nabla u,\nabla e_{\alpha}\rangle+\nabla\cdot\langle \Delta u,\nabla e_{\alpha} \rangle\\
+\sum_{\beta}\Big(\nabla\cdot(\langle\Delta u, e_{\beta}\rangle \langle\nabla e_{\alpha}, e_{\beta}\rangle )-
\langle(\Delta u)^T, e_{\beta}\rangle \langle\nabla e_{\alpha}, e_{\beta}\rangle\Big)\\
-\sum_{\beta}\langle\Delta u, e_{\beta}\rangle \langle \mathbb{B}(u)(e_{\alpha},\nabla u), e_{\beta}\rangle
+\Big\langle\nabla(\mathbb{B}(u)(\nabla u,\nabla u)), \nabla e_{\alpha}\Big\rangle.
\end{cases}
\end{equation}
As in \cite{W2} Lemma 4.3,  we let $\tilde{u}\in W^{2,2}(\mathbb R^4,\mathbb R^K), \{\tilde{e}_{\alpha}\}_{\alpha=1}^n, \tilde{A}$ and $\tilde{h}$ be the extensions of
$u, \{{e}_{\alpha}\}_{\alpha=1}^n, A$ and ${h}$ from $B_\frac12$ to $\mathbb R^4$ such that
$ |\tilde{e}_{\alpha}|\le 1, \tilde{A}^{\alpha\beta}=\langle\nabla \tilde{e}_{\alpha},\tilde{e}_{\beta}\rangle$, and
\begin{equation}\label{extend}
\begin{cases}
\|\nabla^2\widetilde u\|_{L^{2}(\mathbb R^4)}
\le C\|\nabla^2 u\|_{L^{2}(B_\frac12)},\quad
\|\nabla\widetilde u\|_{L^{4}(\mathbb R^4)}\le C\|\nabla u\|_{L^{4}(B_\frac12)};\\
\|\nabla\widetilde u\|_{L^{4,\infty}(\mathbb R^4)}\le C\|\nabla u\|_{L^{4,\infty}(B_\frac12)},
\quad\|\nabla^2\widetilde u\|_{L^{2,\infty}(\mathbb R^4)}\le C\|\nabla^2 u\|_{L^{2,\infty}(B_\frac12)};\\
\|\tilde{A}\|_{L^{4,1}(\mathbb R^4)}+\|\nabla \tilde{A}\|_{L^{2,1}(\mathbb R^4)}
\le C(\|\nabla u \|_{L^{4}(B_1)}+\|\nabla^2 u\|_{L^2(B_1)}); \\
\|\tilde{h}\|_{L^q(\mathbb R^4)}\le C \|h\|_{L^q(B_1)}, \ \ \ \forall 1\le q\le p.
\end{cases}
\end{equation}
Let $\Gamma(x-y)=c_4\ln|x-y|$ be the fundamental solution of $\Delta^2$ in $\mathbb R^4$. Set
\begin{eqnarray}\label{solution representation}
W_{\alpha}(x)&=&\int_{\mathbb R^4}\Gamma(x-y)G(\tilde{u},\tilde{e}_{\alpha})(y)dy+\int_{\mathbb R^4}\Gamma(x-y)\langle\tilde{h},\tilde{e}_{\alpha}\rangle(y)dy\nonumber\\
&=&J^{\alpha}_1(x)+J^{\alpha}_2(x)\  \ x\in\ \ \mathbb R^4.
\end{eqnarray}
Then we have
\begin{equation}\label{W-eqn}
\Delta^2 W_{\alpha}=G(\tilde{u},\tilde{e}_{\alpha})+\langle \tilde{h}, \tilde{e}_{\alpha}\rangle, \ {\rm{in}}\ \mathbb R^4.
\end{equation}
For $J_2^\alpha$, observe that
$$\nabla^4 J_2^\alpha(x)=\int_{\mathbb R^4}\nabla^4_y\Gamma(x-y) \langle\widetilde h, \widetilde{e}_\alpha\rangle(y)\,dy.$$
Hence by Calderon-Zygmund's $W^{4,p}$-theorey, we have $J_2^\alpha\in W^{4,p}(\mathbb R^4)$ and
\begin{equation}\label{J2-estimate}
\Big\|\nabla^4 J_2^\alpha\Big\|_{L^{p}(\mathbb R^4)}\le C\Big\|\widetilde h\Big\|_{L^p(\mathbb R^4)}
\le C\Big\|h\Big\|_{L^p(B_1)}.
\end{equation}
By Sobolev embedding theorem, we have that $\nabla^2 J_1^\alpha\in L^{\bar p}(B_1)$,
with $\bar p=\frac{2p}{2-p}$ for $1<p<2$ or $\bar p$ is any finite number for $p\ge 2$, and
\begin{equation}\label{J2-estimate1}
\Big\|\nabla^2 J_2^\alpha\Big\|_{L^{\bar p}(B_1)}\le C\Big\|h\Big\|_{L^p(B_1)}.
\end{equation} 
By H\"older inequality, (\ref{J2-estimate}) and (\ref{J2-estimate1}) imply that  for any $\theta\in (0,\frac12)$, it holds
\begin{equation}\label{J2-estimate2}
\Big\|\nabla J_2^\alpha\Big\|_{L^{4,\infty}(B_\theta)}
+\Big\|\nabla^2 J_2^\alpha\Big\|_{L^{2,\infty}(B_\theta)}
\le C\theta^{4(1-\frac{1}{p})}\Big\|\nabla^2 J_2^\alpha\Big\|_{L^{\bar p}(B_\theta)}
\le C\theta^{4(1-\frac{1}{p})}\Big\|h\Big\|_{L^p(B_1)}.
\end{equation}
For $J_1^\alpha$, we follow exactly \cite{W2}  Lemma 4.4, Lemma 4.5, and Lemma 4.6 to obtain that
$\nabla^2 J_1^\alpha\in L^{2,\infty}(\mathbb R^4)$ and
\begin{equation}\label{J1-estimate}
\Big\|\nabla^2 J_1^\alpha\Big\|_{L^{2,\infty}(\mathbb R^4)}
\le C\varepsilon_0\Big(\Big\|\nabla u\Big\|_{L^{4,\infty}(B_\frac12)}+\Big\|\nabla^2 u\Big\|_{L^{2,\infty}(B_\frac12)}\Big).
\end{equation}
By Sobolev embedding theorem in Lorentz spaces, (\ref{J1-estimate}) yields
\begin{equation}\label{J1-estimate1}
\Big\|\nabla J_1^\alpha\Big\|_{L^{4,\infty}(\mathbb R^4)}\le C\Big\|\nabla ^2J_1^\alpha\Big\|_{L^{2,\infty}(\mathbb R^4)}.
\end{equation}
Combining (\ref{J1-estimate}) with (\ref{J1-estimate1}) yields
\begin{equation}\label{J1-estimate2}
\Big\|\nabla J_1^\alpha\Big\|_{L^{4,\infty}(\mathbb R^4)}+\Big\|\nabla^2 J_1^\alpha\Big\|_{L^{2,\infty}(\mathbb R^4)}
\le C\varepsilon_0\Big(\Big\|\nabla u\Big\|_{L^{4,\infty}(B_\frac12)}+\Big\|\nabla^2 u\Big\|_{L^{2,\infty}(B_\frac12)}\Big).
\end{equation}

Now, as in \cite{W2} Lemma 4.6,  we consider the Hodge decomposition of the 1-form $\langle d\tilde{u},\tilde{e}_{\alpha} \rangle$. It is well-known  \cite{IM} that there exist $F_{\alpha}\in W^{1,4}(\mathbb R^4)$ and 
$H_{\alpha}\in W^{1,4}(\mathbb R^4, \wedge^2\mathbb R^4)$ such that 
\begin{equation}\label{hodge}
\langle d\tilde{u},\tilde{e}_{\alpha} \rangle = dF_{\alpha}+d^*H_{\alpha},\quad dH_{\alpha}= 0\ \ \ \ {\rm{in}}\ \ \ \mathbb R^4,
\end{equation}
\begin{equation}\label{4-4-weak-estimate}
\Big\|\nabla F_{\alpha}\Big\|_{L^{4}(\mathbb R^4)}+\Big\|\nabla H_{\alpha}\Big\|_{L^{4}(\mathbb R^4)} \le
C\Big\|\nabla \tilde{u}\Big\|_{L^{4}(\mathbb R^4)}\le C\Big\|\nabla u\Big\|_{L^4(B_\frac12)}.
\end{equation}
It is easy to see that $H_{\alpha}$ satisfies
\begin{equation}\label{H-eqn}
\Delta H_{\alpha}= d\widetilde{u}\wedge d\widetilde{e}_{\alpha} \qquad {\rm{in}}\ \ \ \ \mathbb R^4,
\end{equation}
and
\begin{equation}\label{F_eqn}
\Delta^2 F_{\alpha}=\Delta\nabla\cdot\langle\nabla u,e_{\alpha} \rangle=\Delta^2 W_{\alpha}
\ \ \ {\rm{in}}\ \ \ B_\frac12.
\end{equation}
By Calderon-Zygmund's $L^{r,s}$-theory and Sobolev embedding theorem,
 we have that $\nabla^2 H_\alpha\in L^{2,\infty}(\mathbb R^4)$ and
\begin{equation}\label{H-estimate}
\Big\|\nabla H_\alpha\Big\|_{L^{4,\infty}(\mathbb R^4)}+\Big \|\nabla^2 H_{\alpha}\Big \|_{L^{2,\infty}(\mathbb R^4)}
 \le  C\Big\|d\widetilde{u}\wedge d\widetilde{e}_{\alpha}\Big\|_{L^{2,\infty}(\mathbb R^4)}\le C\varepsilon_0\Big\|\nabla u\Big\|_{L^{4,\infty}(B_\frac12)}.
\end{equation}
By (\ref{F_eqn}), we have that $F_\alpha-W_\alpha$ is a biharmonic function on $B_\frac12$. By the standard estimate of biharmonic
functions, we have (see \cite{W2} Lemma 4.7) that for any $\theta\in (0,\frac12)$, it holds
\begin{eqnarray}\label{F-W-estimate}
&&\Big\|\nabla(F_\alpha-W_\alpha)\Big\|_{L^{4,\infty}(B_\theta)}
+\Big\|\nabla^2(F_\alpha-W_\alpha)\Big\|_{L^{2,\infty}(B_\theta)}\nonumber\\
&\le& C\theta\Big(\Big\|\nabla(F_\alpha-W_\alpha)\Big\|_{L^{4,\infty}(B_1)}
+\Big\|\nabla^2(F_\alpha-W_\alpha)\Big\|_{L^{2,\infty}(B_1)}\Big).
\end{eqnarray}
Putting (\ref{J2-estimate2}), (\ref{J1-estimate2}), (\ref{H-estimate}), and (\ref{F-W-estimate}) together, we can argue, similar to
\cite {W2} page 84, to reach that for any $\theta\in (0,\frac12)$, it holds
\begin{eqnarray}\label{u22-estimate}
\Big\|\nabla u\Big\|_{L^{4,\infty}(B_\theta)}
+\Big\|\nabla^2 u\Big\|_{L^{2,\infty}(B_\theta)}
&\le& C(\varepsilon+\theta)
\Big(\Big\|\nabla u\Big\|_{L^{4,\infty}(B_1)}
+\Big\|\nabla^2u\Big\|_{L^{2,\infty}(B_1)}\Big)\nonumber\\
&&+C\theta^{4(1-\frac{1}{p})}\Big\|h\Big\|_{L^p(B_1)}.
\end{eqnarray}
It is readily seen that for any $\alpha\in (0, \min\{1, 4(1-\frac1{p})\})$,
we can choose both $\theta\in (0,\frac12)$ and $\varepsilon_0\in (0,1)$ sufficiently small so that
\begin{equation}\label{u22-estimate1}
\Big\|\nabla u\Big\|_{L^{4,\infty}(B_\theta)}
+\Big\|\nabla^2 u\Big\|_{L^{2,\infty}(B_\theta)}
\le \theta^\alpha\Big(\Big\|\nabla u\Big\|_{L^{4,\infty}(B_1)}
+\Big\|\nabla^2 u\Big\|_{L^{2,\infty}(B_1)}+\Big\|h\Big\|_{L^p(B_1)}\Big).
\end{equation}
In fact, by iterating (\ref{u22-estimate1}) finitely many times  on $B_\frac12(x)$ for $x\in B_\frac12$, we would have
that for $0<r\le \frac14$, 
\begin{equation}\label{u22-estimate2}
\Big\|\nabla u\Big\|_{L^{4,\infty}(B_r(x))}
+\Big\|\nabla^2 u\Big\|_{L^{2,\infty}(B_r(x))}
\le Cr^\alpha \Big(\Big\|\nabla u\Big\|_{L^{4,\infty}(B_1)}
+\Big\|\nabla^2 u\Big\|_{L^{2,\infty}(B_1)}+\Big\|h\Big\|_{L^p(B_1)}\Big).
\end{equation}
Since $L^q(B_r(x))\subset L^{s,\infty}(B_r(x))$ for any $1\le q<s$, (\ref{u22-estimate2}) implies that for any
$1<q<2$, 
\begin{equation} \label{u22-estimate3}
r^{2q-4}\int_{B_r(x)}(|\nabla u|^{2q}+|\nabla^2 u|^q)
\le Cr^{2\alpha} \Big(\Big\|\nabla u\Big\|_{L^{4}(B_1)}^4
+\Big\|\nabla^2 u\Big\|_{L^{2}(B_1)}^2+\Big\|h\Big\|_{L^p(B_1)}^p\Big).
\end{equation}
This, with the help of Morrey's decay lemma, immediately implies  that $u\in C^\alpha(B_\frac12)$ and
(\ref{holder_est}) holds.
This completes the proof of Lemma 2.1. \qed\\

In order to show $\varepsilon_0$-compactness and Pohozaev argument for approximate biharmonic maps, we will establish the
higher order Sobolev type regularity for approximate biharmonic maps.

The proof utilizes Adams' Reisz potential estimate between Morry spaces, we briefly recall Morrey spaces and Adams' estimates (see
\cite{A} and \cite{W3} for more details).
For an open set $U\subset\mathbb R^4$, $1\le p<+\infty$, $0<\lambda\le 4$, the Morrey space $M^{p,\lambda}(U)$ is defined by
\begin{equation}\label{Morrey}
M^{p,\lambda}(U)
=\Big\{f\in L^p(U):
\|f\|_{M^{p,\lambda}}^p
=\sup_{B_r\subset U} r^{\lambda-4}\int_{B_r}|f|^p<+\infty \Big\}.
\end{equation}
The weak Morrey space $M^{p,\lambda}_*(U)$ is  the
set of functions $f\in L^{p,\infty}(U)$ satisfying
\begin{equation}\label{weak Morrey}
\|f\|^p_{M^{p,\lambda}_*(U)}\equiv \sup_{B_r\subset U}\Big\{\rho^{\lambda-n}\|f\|^p_{L^{p,\infty}(B_r)}\Big\}<\infty.
\end{equation}

For $ 0<\beta<4$, let $I_\beta(f)$ be the Riesz potential of order  defined by
\begin{equation} \label{reisz}
I_{\beta}(f)(x)\equiv\int_{\mathbb{R}^4}\frac
{f(y)}{|x-y|^{4-\beta}} \,dy, \qquad x\in \mathbb{R}^4.
\end{equation}
Recall Adams' estimate \cite{A} in dimension $4$:
\begin{lemma}
(1)\  For any $\beta>0, 0 <\lambda\le 4, 1 < p <\frac
{\lambda}{\beta}$, if $f\in M^{p,\lambda}(\mathbb{R}^4)$, then $I_{\beta}(f)\in
M^{\tilde{p},\lambda}(\mathbb{R}^4)$, where $\tilde{p}=\frac
{\lambda p}{\lambda-p\beta} $. Moreover,
\begin{equation}\label{Riesz-Morrey}
\|I_{\beta}(f)\|_{M^{\tilde{p},\lambda}(\mathbb{R}^4)}\le
C\|f\|_{M^{p,\lambda}(\mathbb{R}^4)}.
\end{equation}
(2)\  For $0<\beta<\lambda\le 4$, if $f\in M^{1,\lambda}(\mathbb{R}^4)$, then $I_{\beta}(f)\in
M_*^{\frac{\lambda}{\lambda-\beta},\lambda}(\mathbb{R}^4)$. Moreover,
\begin{equation}\label{Riesz-weak-Morrey}
\|I_{\beta}(f)\|_{M_*^{\frac{\lambda}{\lambda-\beta},\lambda}(\mathbb{R}^4)}\le
C\|f\|_{M^{1,\lambda}(\mathbb{R}^4)}.
\end{equation}
\end{lemma}

Now we are ready to prove
\begin{lemma} \label{higher integrality} There exists $\varepsilon_0>0$ such that  if
$u\in W^{2,2}(B_1,N)$ is an approximate biharmonic map with its bi-tension field $h\in L^p(B_1)$ for $p>1$, and
satisfies
\begin{equation}\label{small_energy1}
\int_{B_1}|\nabla u|^4+|\nabla^2 u|^2\le\varepsilon_0^2.
\end{equation}
Then $u\in W^{4,p}(B_{\frac 18},N)$ and
\begin{equation}\label{4-p-integ}
\Big\|\nabla^4u\Big\|_{L^{p}(B_{\frac 18})}\le C\Big(\Big\|\nabla u\Big\|_{L^4(B_1)}+\Big\|\nabla^2
u\Big\|_{L^2(B_{1})}+\Big\|h\Big\|_{L^p(B_1)}\Big).
\end{equation}
\end{lemma}
\pf By (\ref{holder_est}) of Lemma 2.1, we have that for any $\alpha\in (0,\min\{1, 4(1-\frac{1}{p})\})$,
\begin{equation}\label{oscillation}
\mbox{osc}_{B_r(x)} u\in Cr^\alpha, \ \forall B_r(x)\subset B_\frac12.
\end{equation}
We now divide the proof into three steps. \\ 
\noindent{\bf Step1}. There exists $\alpha_0 \in (0, \alpha]$ such that
$\nabla u\in M^{4,4-2\alpha_0}(B_\frac12), \ \nabla^2 u\in M^{2,4-2\alpha_0}(B_\frac12)$, and
\begin{equation}\label{morrey0}
\Big\|\nabla u\Big\|_{M^{4,4-2\alpha_0}(B_\frac12)}+\Big\|\nabla^2 u\Big\|_{M^{2,4-2\alpha_0}(B_\frac12)}
\le C\Big(\Big\|\nabla u\Big\|_{L^4(B_1)}+\Big\|\nabla^2
u\Big\|_{L^2(B_{1})}+\Big\|h\Big\|_{L^p(B_1)}\Big).
\end{equation}
Observe that (\ref{morrey0}) is a refined version of (\ref{u22-estimate2}) and (\ref{u22-estimate3}). It is obtained by the hole filling
argument as follows. For any $B_r(x)\subset B_\frac12$, let $\phi\in C_0^\infty(B_r(x))$ such that
$0\le \phi\le 1$, $\phi\equiv 1$ on $B_{\frac{r}2}(x)$. Multiplying (\ref{approx_biharm}) by $\phi(u-u_{x,r})$, where
$u_{x,r}$ is the average of $u$ on $B_r(x)$, and integrating over $B_r(x)$, we would have
\begin{eqnarray}\label{hole1}
&&\int_{B_r(x)}|\Delta(\phi(u-u_{x,r}))|^2\nonumber\\
&\le& \int_{B_r(x)} \Delta((1-\phi)(u-u_{x,r}))\cdot\Delta(\phi(u-u_{x,r}))
+C\int_{B_r(x)}|\nabla u|^2 |\Delta(\phi(u-u_{x,r}))|\nonumber\\
&+&C\int_{B_r(x)}|\nabla^2 u||\nabla u||\nabla(\phi(u-u_{x,r})|+C\Big(\int_{B_r(x)}|\Delta u|^2+|h|\Big)\mbox{osc}_{B_r(x)} u.
\end{eqnarray}
It is not hard to see that by H\"older inequality, Sobolev's inequality and (\ref{oscillation}), (\ref{hole1}) implies
\begin{equation}\label{hole2}
\int_{B_{\frac{r}2}(x)}(|\nabla u|^4+|\nabla^ 2u|^2)
\le \theta\int_{B_r(x)}(|\nabla u|^4+|\nabla^2 u|^2)+Cr^\alpha,
\end{equation}
where $\theta=\frac{C}{C+1} <1$. Now we can iterate (\ref{hole2}) finitely many times and achieve that
there exists $\alpha_0\in (0,\alpha)$ such that
\begin{equation}
\sup_{B_r(x)\subset B_\frac12} r^{-2\alpha_0}\int_{B_r(x)}(|\nabla u|^4+|\nabla^2 u|^2)
\le C\Big(\|\nabla u\|_{L^4(B_1)}^4+\|\nabla^2 u\|_{L^2(B_1)}^2+\|h\|_{L^p(B_1)}^p\Big).
\end{equation}
This yields (\ref{morrey0}).\\
\noindent{\bf Step 2}.  Set $\bar p>2$ by
$$
\bar p=\begin{cases}\frac{2p}{2-p} & {\rm{if}} \ 1<p<2\\
{\rm{any}}\  2<q<+\infty & {\rm{if}} \ p\ge 2.
\end{cases}
$$
Then $u\in W^{2,\bar p}(B_\frac14)$ and
\begin{equation}\label{bar-p}
\Big\|\nabla u\Big\|_{L^{2\bar p}(B_\frac14)}+\Big\|\nabla^2 u\Big\|_{L^{\bar p}(B_\frac14)}\le C\Big(\Big\|\nabla u\Big\|_{L^4(B_1)}
+\Big\|\nabla^2 u\Big\|_{L^2(B_1)}+\Big\|h\Big\|_{L^p(B_1)}\Big).
\end{equation}
To show (\ref{bar-p}),  let $\widetilde u,\widetilde h:\mathbb R^4\to \mathbb R^{K}$ be extensions
of $u$ and $h$ on $B_\frac12$ such that
\begin{equation}\label{u_extension}
\|\nabla^2 \widetilde u\|_{M^{2,4-2\alpha_0}(\mathbb R^4)}
\le C\|\nabla^2 u\|_{M^{2,4-2\alpha_0}(B_{\frac12})},\quad \|\nabla \widetilde u\|_{M^{4,4-2\alpha_0}(\mathbb R^4)}
\le C\|\nabla u\|_{M^{4,4-2\alpha_0}(B_{\frac 12})},
\end{equation}
and
\begin{equation}\label{h_extension}
\|\widetilde{h}\|_{L^p(\mathbb R^4)}\le C\|h\|_{L^{p}(B_1)}.
\end{equation}
Define $w:\mathbb R^4\to\mathbb R^K$ by 
\begin{eqnarray}\label{w-representation}
w(x)&=&\int_{\mathbb R^4}\Gamma(x-y)\widetilde{h}(y)\,dy
+\int_{\mathbb R^4}\Delta_y \Gamma(x-y) (\mathbb B(\tilde{u})(\nabla \tilde{u},\nabla \tilde{u}))(y)\, dy\nonumber\\
&-&2\int_{\mathbb R^4}\nabla_y \Gamma(x-y)\langle\Delta \tilde{u},
\nabla(\mathbb P(\tilde{u}))\rangle(y)\, dy
-\int_{\mathbb R^4}\Gamma(x-y)\langle\Delta(\mathbb P(\tilde{u})),\Delta \tilde{u}\rangle(y)\,dy\nonumber\\
&=& w_1(x)+w_2(x)+w_3(x)+w_4(x),  \ x\in\mathbb R^4.
\end{eqnarray}
Then it is readily seen that 
\begin{equation}\label{biharm_fun}
\Delta^2 (u-w)=0  \ \ \ {\rm{on}}\ \ \ B_{\frac 12},
\end{equation}
or $u-w$ is a biharmonic function on $B_\frac12$.

Now we estimate $w_i$, $1\le i\le 4$, as follows. For $w_1$, by Calderon-Zygmund's $W^{4,p}$-theory
we have that $w_1\in W^{4,p}(\mathbb R^4)$ so that $\nabla^2 w_1\in L^{\bar p}(B_\frac12)$ and
\begin{equation}\label{w1-est}
\Big\|\nabla^2 w_1\Big\|_{L^{\bar p}(B_\frac12)}\le C\|h\|_{L^p(B_1)}.
\end{equation}
For $w_3$, since $\displaystyle|\nabla^2 \widetilde u||\nabla\widetilde u|\in M^{\frac43, 4-2\alpha_0}(\mathbb R^4)$, 
$|\nabla w_3|\le CI_2(|\nabla^2 \widetilde u||\nabla\widetilde u|)$ and
$\displaystyle
|\nabla^2w_3|\le CI_1(|\nabla^2 \widetilde u||\nabla\widetilde u|),$
Lemma 2.2 implies that $\nabla w_3\in M^{\frac{4(2-\alpha_0)}{2-3\alpha_0}, 4-2\alpha_0}(\mathbb R^4)$,
$\displaystyle\nabla^2 w_3\in M^{\frac{2(4-2\alpha_0)}{4-3\alpha_0}, 4-2\alpha_0}(\mathbb R^4)$,  and
\begin{eqnarray}\label{w3-est}
&&\Big\|\nabla w_3\Big\|_{M^{\frac{4(2-\alpha_0)}{2-3\alpha_0}, 4-2\alpha_0}(\mathbb R^4)}
+\Big\|\nabla^2 w_3\Big\|_{M^{\frac{2(4-2\alpha_0)}{4-3\alpha_0}, 4-2\alpha_0}(\mathbb R^4)}\nonumber\\
&\le& C\Big\||\nabla^2 \widetilde u||\nabla\widetilde u|\Big\|_{M^{\frac43, 4-2\alpha_0}(\mathbb R^4)}
\le C\Big\|\nabla^2 \widetilde u\Big\|_{M^{2, 4-2\alpha_0}(\mathbb R^4)}
\Big\|\nabla\widetilde u\Big\|_{M^{4, 4-2\alpha_0}(\mathbb R^4)}\nonumber\\
&\le& C\Big\|\nabla^2 u\Big\|_{M^{2, 4-2\alpha_0}(B_\frac12)}
\Big\|\nabla u\Big\|_{M^{4, 4-2\alpha_0}(B_\frac12)}\nonumber\\
&\le& C\Big(\Big\|\nabla u\Big\|_{L^4(B_1)}+\Big\|\nabla^2
u\Big\|_{L^2(B_{1})}+\Big\|h\Big\|_{L^p(B_1)}\Big).
\end{eqnarray}
For $w_4$, it is easy to see that $\displaystyle|\nabla w_4|\le CI_3(|\nabla^2 \widetilde u|^2+|\nabla\widetilde u|^4)$
and $\displaystyle|\nabla^2 w_4|\le CI_2(|\nabla^2 \widetilde u|^2+|\nabla\widetilde u|^4)$. Since
$\displaystyle(|\nabla^2\widetilde u|^2+|\nabla\widetilde u|^4)\in M^{1,4-2\alpha_0}(\mathbb R^4)$, Lemma 2.2 implies
that $\displaystyle|\nabla w_4|\in M_*^{\frac{4-2\alpha_0}{1-2\alpha_0}, 4-2\alpha_0}(\mathbb R^4)$
and $\displaystyle|\nabla^2 w_4|\in M_*^{\frac{2-\alpha_0}{1-\alpha_0}, 4-2\alpha_0}(\mathbb R^4)$ and
\begin{eqnarray}\label{w4-est}
&&\Big\|\nabla w_4\Big\|_{ M_*^{\frac{4-2\alpha_0}{1-2\alpha_0}, 4-2\alpha_0}(\mathbb R^4)}
+\Big\|\nabla^2 w_4\Big\|_{ M_*^{\frac{2-\alpha_0}{1-\alpha_0}, 4-2\alpha_0}(\mathbb R^4)}\nonumber\\
&\le& C\Big(\Big\|\nabla u\Big\|_{L^4(B_1)}+\Big\|\nabla^2
u\Big\|_{L^2(B_{1})}+\Big\|h\Big\|_{L^p(B_1)}\Big).
\end{eqnarray}
For $w_2$, since $\displaystyle|\nabla w_2|\le CI_1(|\nabla \widetilde u|^2)$,  $|\nabla\widetilde u|^2\in
M^{2,4-2\alpha_0}(\mathbb R^4)$, Lemma 2.2 implies that
$\displaystyle |\nabla w_2|\in M^{\frac{4-2\alpha_0}{1-\alpha_0}, 4-2\alpha_0}(\mathbb R^4)$ and
\begin{eqnarray}\label{w2-est}
\Big\|\nabla w_2\Big\|_{M^{\frac{4-2\alpha_0}{1-\alpha_0}, 4-2\alpha_0}(\mathbb R^4)}
&\le& C\Big\||\nabla\widetilde u|^2\Big\|_{M^{2,4-2\alpha_0}(\mathbb R^4)}\nonumber\\
&\le& C\Big(\Big\|\nabla u\Big\|_{L^4(B_1)}+\Big\|\nabla^2
u\Big\|_{L^2(B_{1})}+\Big\|h\Big\|_{L^p(B_1)}\Big).
\end{eqnarray}
It is not hard to see from (\ref{biharm_fun}) and the standard estimate on biharmonic function,
and the estimates (\ref{w1-est}),  (\ref{w3-est}), (\ref{w4-est}), and (\ref{w2-est}) that
there exist $1<q<\min\Big\{\frac{p}{4-3p}, \frac{4-2\alpha_0}{4(1-\alpha_0)}, \frac{2-\alpha_0}{2-3\alpha_0}\Big\}$ and
$0<\alpha_1\le\min\Big\{\alpha_0, \frac{(4-3p)q}{2p}\Big\}$ such that 
$\nabla u\in M^{4q,4-2\alpha_1}(B_\frac38)$ and $\nabla^2 u\in M^{2q, 4-2\alpha_1}(B_\frac38)$, and
\begin{equation}\label{2q-est1}
\Big\|\nabla u\Big\|_{M^{4q,4-2\alpha_1}(B_\frac38)}
+\Big\|\nabla^2 u\Big\|_{M^{2q, 4-2\alpha_1}(B_\frac38)}
\le C\Big(\Big\|\nabla u\Big\|_{L^4(B_1)}+\Big\|\nabla^2
u\Big\|_{L^2(B_{1})}+\Big\|h\Big\|_{L^p(B_1)}\Big).
\end{equation}
With (\ref{2q-est1}), we can repeat the same argument to bootstrap the integrablity of $\nabla^2 u$ and finally get that
$\nabla^2 u\in L^{\bar p}(B_\frac14)$ and (\ref{bar-p}) holds. \\
\noindent{\bf Step 3}. $\nabla^4 u\in L^p(B_\frac18)$ and
\begin{equation}\label{4p-est}
\Big\|\nabla^4 u\Big\|_{L^p(B_\frac18)}
\le C\Big(\Big\|\nabla u\Big\|_{L^4(B_1)}+\Big\|\nabla^2
u\Big\|_{L^2(B_{1})}+\Big\|h\Big\|_{L^p(B_1)}\Big).
\end{equation}
To prove (\ref{4p-est}), first observe that the equation (\ref{approx_biharm}) can be written as
\begin{equation}\label{approx_biharm0}
\Delta^2 u=\mbox{div} (E(u))+G(u)+h,\ \  \  {\rm{in}}\ \ \ B_\frac14,
\end{equation}
where $E(u)=\nabla(B(u)(\nabla u,\nabla u))+2\langle\Delta u, \nabla(\mathbb P(u))\rangle$ and $G(u)=-\langle \Delta(\mathbb P(u)),\Delta u\rangle$ so
that 
$$|E(u)|\le C(|\nabla u|^3+|\nabla^2 u||\nabla u|), 
\ \ |G(u)|\le C(|\nabla^2u|^2+|\nabla u|^4).$$
By (\ref{bar-p}) and Sobolev's inequality, we have $\displaystyle\nabla u\in L^{\frac{4p}{4-3p}}(B_\frac14)$ so that
$\displaystyle G(u)\in L^{\frac{p}{2-p}}(B_\frac14)$ and $\displaystyle E(u)\in L^{\frac{4p}{8-5p}}(B_\frac14)$. 
Note that we can write $u=u_1+u_2+u_3+u_4$ in $B_\frac14$, where
\begin{equation}\label{u_1-decomposition}
\Delta^2 u_1=G(u) \ {\rm{in}}\ B_\frac14; \ (u_1,\nabla u_1)=(0,0) \ {\rm{on}}\ \partial B_\frac14,
\end{equation}
\begin{equation}\label{u_2-decomposition}
\Delta^2 u_2= h \ {\rm{in}}\ B_\frac14; \ (u_2,\nabla u_2)=(0,0) \ {\rm{on}}\ \partial B_\frac14,
\end{equation}
\begin{equation}\label{u_3-decomposition}
\Delta^2 u_3= \hbox{div}(E(u)) \ B_\frac14; \ (u_3,\nabla u_3)=(0,0) \ {\rm{on}}\ \partial B_\frac14,
\end{equation}
and
\begin{equation}\label{u_4-decomposition}
\Delta^2 u_4= 0 \ B_\frac14; \ (u_3,\nabla u_3)=(u, \nabla u) \ {\rm{on}}\ \partial B_\frac14.
\end{equation}
By Calderon-Zygmund's $L^q$-theory, we have 
\begin{equation}\label{u1-u2-estimate}
\Big\|u_1\Big\|_{W^{4, \frac{p}{2-p}}(B_\frac14)}
+\Big\|u_2\Big\|_{W^{4,p}(B_\frac14)}\le C(\Big(\Big\|\nabla u\Big\|_{L^4(B_1)}+\Big\|\nabla^2
u\Big\|_{L^2(B_{1})}+\Big\|h\Big\|_{L^p(B_1)}\Big),
\end{equation}
and
\begin{equation}\label{u3-u4-estimate}
\Big\|u_3\Big\|_{W^{3, \frac{4p}{4-p}}(B_\frac14)}
+\Big\|u_4\Big\|_{W^{4,p}(B_\frac15)}\le C(\Big(\Big\|\nabla u\Big\|_{L^4(B_1)}+\Big\|\nabla^2
u\Big\|_{L^2(B_{1})}+\Big\|h\Big\|_{L^p(B_1)}\Big),
\end{equation}
In particular, we can conclude that $u\in W^{3, \frac{4p}{4-p}}(B_\frac15)$ and
\begin{equation}\label{3q-est}
\Big\|u\Big\|_{W^{3, \frac{4p}{4-p}}(B_\frac15)}\le C(\Big(\Big\|\nabla u\Big\|_{L^4(B_1)}+\Big\|\nabla^2
u\Big\|_{L^2(B_{1})}+\Big\|h\Big\|_{L^p(B_1)}\Big).
\end{equation}
By H\"older inequality, (\ref{3q-est}) then implies 
$$
|\hbox{div}(E(u))|\le C(|\nabla^3 u||\nabla u|+|\nabla^2 u|^2+|\nabla u|^3+|\nabla u|^4)\in L^{\frac{p}{2-p}}(B_\frac15).$$
Hence applying $W^{4,q}$-estimate of (\ref{u_3-decomposition}) yields that
$u_3\in W^{4,\frac{p}{2-p}}(B_\frac18)$ and 
\begin{equation}\label{u3-estimate}
\Big\|u_3\Big\|_{W^{4, \frac{p}{2-p}}(B_\frac18)}\le C(\Big(\Big\|\nabla u\Big\|_{L^4(B_1)}+\Big\|\nabla^2
u\Big\|_{L^2(B_{1})}+\Big\|h\Big\|_{L^p(B_1)}\Big).
\end{equation}
Since $\frac{p}{2-p}>p$, by combining (\ref{u3-estimate}) with (\ref{u1-u2-estimate}) and (\ref{u3-u4-estimate}), we finally obtain that
$u\in W^{4,p}(B_\frac18)$ and (\ref{4p-est}) holds.
This completes the proof of Lemma 2.3.   \qed

\section {Blow up analysis and energy inequality}
\setcounter{equation}{0}
\setcounter{theorem}{0}

This section is devoted to $\epsilon_0$-compactness lemma and
preliminary steps on the blow up analysis of approximate
biharmonic maps with bi-tension fields bounded in $L^p$ for $p>1$.

First we have
\begin{lemma}\label{e-strong} For $n=4$,
there exists an $\epsilon_0>0$ such that if $\{u_k\}\subset
W^{2,2}(B_1, N)$ is a sequence of approximate biharmonic
maps satisfying
\begin{equation}\label{smallness2}
\sup_{k}\Big(\|\nabla u_k\|_{L^4(B_1)}^4+\left\|\nabla^2 u_k\right\|_{L^2(B_1)}^2\Big)\le\epsilon_0,
\end{equation}
and $u_k\rightharpoonup u$ in $W^{2,2}(B_1)$ and
$h_k\rightharpoonup h$ in $L^p(B_1)$ for some $p>1$.
Then $u\in C^0\cap W^{4,p}(\Omega, N)$ is an approximate biharmonic map with bi-tension field
$h$, and
\begin{equation}\label{strong}
\lim_{k\rightarrow\infty}
\Big\|u_k-u\Big\|_{W^{2,2}(B_\frac12)}=0.
\end{equation}
\end{lemma}
\pf The first assertion follows easily from (\ref{approx_biharm}) and (\ref{strong}). To show (\ref{strong}), it suffices
to show that $\{u_k\}$ is a Cauchy sequence in $W^{2,2}(B_\frac12)$.
By (\ref{smallness2}) and Lemma \ref{holder}, there exist
$\alpha\in (0,1)$ and $q>2$ such that
$$\sup_k\Big[ \Big\|u_k\Big\|_{C^\alpha(B_\frac34)}
+\Big\|\nabla^2 u_k\Big\|_{L^q(B_\frac34)}\Big]\le C.$$
Hence we may assume that
$$\lim_{k, l\rightarrow\infty} \Big\|u_k-u_l\Big\|_{L^\infty(B_\frac34)}=0.$$
For $\eta\in C_0^\infty(B_\frac34)$ be a cut-off function of
$B_\frac12$, multiplying the equations of $u_k$ and $u_l$ by
$(u_k-u_l)\phi^2 $ and integrating over $B_1$, we obtain
\begin{eqnarray*}
&&\int_{B_1}|\Delta(u_k-u_l)|^2\phi^2\\
&\le &\int_{B_1}|\Delta(u_k-u_l)|(2|\nabla(u_k-u_l)||\nabla\phi^2|+|u_k-u_l||\Delta\phi^2|)
+\int_{B_1}|h_k-h_l||u_k-u_l|\phi^2\\
&&+3\int_{B_1}(|\Delta u_l|^2 +|\Delta u_k|^2)|u_k-u_l|\phi^2\\
&&+4\int_{B_1}|\nabla^2 u_k||\nabla u_k||\nabla(u_k(u_k-u_l)\phi^2)|\\
&&+4\int_{B_1}|\nabla^2 u_l||\nabla u_l||\nabla(u_l(u_k-u_l)\phi^2)|\\
&=&I+II+III+IV+V.
\end{eqnarray*}
It is easy to see
$$|I|\le C(\|\nabla (u_k-u_l)\|_{L^2(B_\frac34)}
+\|u_k-u_l\|_{L^\infty(B_\frac34)})\rightarrow 0,$$
$$|II|\le C\|h_k-h_l\|_{L^1(B_\frac 34)}\|u_k-u_l\|_{L^\infty(B_\frac34)}\rightarrow 0,$$
$$|III|\le C(\|\nabla^2 u_k\|_{L^2(B_\frac34)}^2
+\|\nabla^2 u_l\|_{L^2(B_\frac34)}^2)\|u_k-u_l\|_{L^\infty(B_\frac34)}\rightarrow 0.$$
For $IV$, observe that for $1<r<4$ with $\frac14+\frac1{q}
+\frac1{r}=1$, we have
\begin{eqnarray*}
|IV|
&\le& C\Big(\|\nabla^2 u_k\|_{L^2(B_\frac34)}\|\nabla u\|_{L^4(B_\frac34)}^2\|u_k-u_l\|_{L^\infty(B_\frac34)}\\
&&+\|\nabla^2 u_k\|_{L^q(B_\frac34)}\|\nabla u_k\|_{L^4(B_\frac34)}
\|\nabla(u_k-u_l)\|_{L^r(B_\frac34)}\Big)
\rightarrow 0,
\end{eqnarray*}
since $\displaystyle\|\nabla(u_k-u_l)\|_{L^r(B_\frac34)}\rightarrow 0$.
Similarly, we can show
$$|V|\rightarrow 0.$$
Hence $\{u_k\}$ is a Cauchy sequence in $W^{2,2}(B_\frac12)$.
This completes the proof.
\qed
\begin{lemma} Under the same assumptions as Theorem \ref{energy_identity}, there exists a finite subset $\Sigma\subset
\Omega$ such that
$u_k\rightarrow u$ in $W^{2,2}_{\rm{loc}}\cap C^0_{\rm{loc}}
(\Omega\setminus\Sigma,N)$. Moreover,
$u\in W^{4,p}\cap C^0(\Omega,N)$ is an approximate
biharmonic map with bi-tension field $h$.
\end{lemma}
\pf Let $\epsilon_0>0$ be given by Lemma \ref{holder}, and define
\begin{equation}\label{concentration}
\Sigma:=\bigcap_{r>0}\Big\{x\in\Omega:
\liminf_{k\rightarrow\infty}\int_{B_r(x)}(|\nabla^2 u_k|^2+|\nabla u_k|^4)>\epsilon_0^2\Big\}.
\end{equation}
Then by a simple covering argument we have that
$\Sigma$ is a finite set and 
$$H^0(\Sigma)\le\frac{1}{\epsilon_0^2}\sup_{k}\int_\Omega
(|\nabla^2 u_k|^2+|\nabla u_k|^4)<+\infty.$$
For any $x_0\in\Omega\setminus\Sigma$, there exists $r_0>0$ such
that
$$\liminf_{k\rightarrow\infty}\int_{B_{r_0}(x_0)}
(|\nabla^2 u_k|^2+|\nabla u_k|^4)\le\epsilon_0^2.$$
Hence Lemma \ref{holder} and Lemma 3.1 imply that there exists
$\alpha\in (0,1)$ such that
$$\Big\|u_k\Big\|_{C^\alpha(B_{\frac{r_0}2}(x_0))}
\le C,$$
so that $u_k\rightarrow u$ in $C^0\cap W^{2,2}(B_{\frac{r_0}2}(x_0))$. This proves that $u_k\rightarrow u$ in
$W^{2,2}_{\rm{loc}}\cap C^0_{\rm{loc}}(\Omega\setminus\Sigma)$.
It is clear that $u\in W^{2,2}(\Omega)$ is an approximate biharmonic map with bi-tension field $h\in L^p(\Omega)$. Applying Lemma \ref{holder} and Lemma 2.3 again, we conclude that $u\in C^0(\Omega,N)\cap W^{4,p}(\Omega, N)$. \qed\\

\noindent{\bf Proof of Theorem\ \ref{energy_identity}}:

\smallskip

The proof of (\ref{energy_id1}) with ``=" replaced by``$\ge$'' is similar to
\cite{W1} Lemma 3.3. Here we sketch it. For any $x_0\in\Sigma$,
there exist $r_0>0$, $x_k\rightarrow x_0$ and $r_k\downarrow 0$ such that
$$
\max_{x\in B_{r_0}(x_0)}
\Big\{\int_{B_{r_k}(x)}(|\nabla^2 u_k|^2+|\nabla u_k|^4)
\Big\}
=\frac{\epsilon_0^2}2=\int_{B_{r_k}(x_k)}(|\nabla^2 u_k|^2+|\nabla u_k|^4).
$$
Define $v_k(x)=u_k(x_k+r_kx): r_k^{-1}\Big(B_{r_0}(x_0)\setminus\{x_k\}\Big)\to N$. 
Then $v_k$ is an approximate biharmonic map, with bi-tension
field $\widetilde{h}_k(x)=r_k^4 h(x_k+r_k x)$, that satisfies
$$\int_{B_1(x)}(|\nabla^2 v_k|^2+|\nabla v_k|^4)\le\frac{\epsilon_0^2}2,
\ \forall x\in r_k^{-1}\Big(B_{r_0}(x_0)\setminus\{x_k\}\Big),
\ {\rm{and}}\ \int_{B_1(0)}(|\nabla^2 v_k|^2+|\nabla v_k|^4)=\frac{\epsilon_0^2}2,$$
and
$$\Big\|\widetilde{h}_k\Big\|_{L^p\Big(r_k^{-1}(B_{r_0}(x_0)\setminus\{x_k\}\Big)}
\le r_k^{4(1-\frac1p)}\Big\|h_k\Big\|_{L^p(\Omega)}
\rightarrow 0.$$
Thus Lemma 2.3 and Lemma 3.1 imply that,
after taking possible subsequences, there exists a nontrivial biharmonic map $\omega:\mathbb R^4\to N$ with
$$\frac{\epsilon_0^2}2\le \int_{\mathbb R^4}(|\nabla^2\omega|^2+|\nabla \omega|^4)<+\infty
$$
such that $v_k\rightarrow \omega$ in
$W^{2,2}_{\rm{loc}}\cap C^0_{\rm{loc}}(\mathbb R^4)$.
Performing such a blow-up argument at all $x_i\in\Sigma$,
$1\le i\le L$, we can find all possible nontrivial biharmonic maps
$\{\omega_{ij}\}\in W^{2,2}(\mathbb R^4)$ for $1\le j\le L_i$,
with $L_i\le CM\epsilon_0^{-2}$. Moreover, by the lower semicontinuity,
we have that the part ``$\ge$" of (\ref{energy_id1}) holds.  

The other half,  ``$\le$'',  of (\ref{energy_id1}) will be proved in the next section.

\section{No hessian energy concentration in the neck region}
\setcounter{equation}{0}
\setcounter{theorem}{0}

In order to show the part ``$\le$" of (\ref{energy_id1}), we need to
show that there is no concentration of hessian energy in the neck region.
This will be done in two steps. The first step is to show that there
is no angular hessian energy concentration in the neck region by comparing
with radial biharmonic functions over dydaic annulus. The
second step is to use the type of almost hessian energy monotonicity inequality,
which is obtained by the Pohozaev type argument, to control the radial component
of hessian energy by the angular component of hessian energy. We require $p>\frac43$ in
both steps.

Suppose that $\{u_k\}\subset W^{2,2}(B_1,N)$ is a sequence of approximate biharmonic maps satisfying for some $p>\frac43$,
\begin{equation}\label{bound}
\int_{B_1}(|\nabla^2 u_k|^2+|\nabla u_k|^4+|h_k|^p)\le C, \ \forall k\ge 1.
\end{equation}
Without loss of generality, we assume that $u_k\rightharpoonup u$ in $W^{2,2}$,   $h_k\rightharpoonup h$ in $L^p$, and 
$ u_k\rightarrow u  $ in $W^{2,2}_{\rm{loc}}(B_1\backslash\{0\},N)$ but not in
$W^{2,2}(B_1, N)$.
Furthermore, as in Ding-Tian \cite{DT}, we may assume that the total number of bubbles generated at $0$ is $L=1$. Then for any $\varepsilon>0$ there exist $r_k\downarrow 0$,
$R\ge 1$ sufficiently large, and $0<\delta\le\varepsilon^{\frac{p}{4(p-1)}}$ so that
for $k$ sufficiently large, the following property holds
\begin{equation}\label{small_energy}
\int_{B_{2\rho}\backslash B_\rho}(|\nabla^2u_k|^2+|\nabla u_k|^4)
\le \epsilon^2,\qquad \forall \  \frac1{16}Rr_k\le \rho \le 16\delta.
\end{equation}

\noindent{\bf Step 1}. {\it Angular hessian energy estimate in the neck region}:
\\

Since $p>\frac43$, it follows from (\ref{small_energy}), Lemma \ref{holder}, Lemma 2.3, and Sobolev embedding theorem
 that for any $\alpha\in \left(0,4(1-\frac{1}{p})\right)$,
$u_k\in C^\alpha\cap W^{4,p}(B_{2\rho}\setminus B_{\rho})$, $\nabla u_k\in C^{0}(B_{2\rho}\setminus B_\rho)$,
and 
\begin{equation}\label{holder annulus}
\Big[u\Big]_{C^{\alpha}(B_{2\rho}\backslash B_{\rho})}
+\Big\|\nabla u_k\Big\|_{L^\infty(B_{2\rho}\setminus B_\rho)}
\le C\Big(\varepsilon+\rho^{4(1-\frac 1p)}\Big)\le C\varepsilon, \ \forall\  \frac{1}8 Rr_k\le \rho\le 8\delta. 
\end{equation}
It follows from Lemma 2.3 that
\begin{equation}\label{343-estimate}
\Big\|\nabla^3 u_k\Big\|_{L^\frac43(B_{2\rho}\setminus B_\rho)}\le C\varepsilon,
\ \forall\  \frac{1}8 Rr_k\le \rho\le 8\delta. 
\end{equation}
To handle the contributions of various boundary terms during the argument, by Fubini's theorem we assume 
that $R>0$ and $\delta>0$ are chosen so that for $k$ sufficiently large, the following property
\begin{equation} \label{Fubini0}
r\int_{\partial B_r}(|\nabla u_k|^4+|\nabla^2 u_k|^2+|\nabla^3 u_k|^\frac43)
\le 8\sup_{k}\int_{B_{2r}\setminus B_{\frac{r}2}}(|\nabla u_k|^4+|\nabla^2 u_k|^2+|\nabla^3 u_k|^\frac43)\le C\varepsilon^2,
\end{equation}
holds for  $r=\frac12 Rr_k, Rr_k,  \delta,  2\delta,  4\delta$. For simplicity, we assume (\ref{Fubini0}) holds for all $k\ge 1$. 
Here we indicate (\ref{Fubini0}) for $r=Rr_k$: set $\widetilde{u}_k=u_k(r_kx): B_{\delta r_k^{-1}}\to N$. Then by 
Fatou's lemma we have
$$\int_{\frac12 R}^{2R}\liminf_{k}\int_{\partial B_r}(|\nabla \widetilde u_k|^4+|\nabla^2 \widetilde u_k|^2+|\nabla^3 \widetilde u_k|^\frac43)
\le \liminf_{k}\int_{B_{2R}\setminus B_{\frac12 R}}(|\nabla \widetilde u_k|^4+|\nabla^2 \widetilde u_k|^2+|\nabla^3 \widetilde u_k|^\frac43).$$
By Fubini's theorem and scalings, this easily imply (\ref{Fubini0}) for $r\approx Rr_k$ (for simplicity, we can assume $r=Rr_k$).

Let $N_k\in \mathbb{N}$ be such that $2^{N_k}=[\frac{2\delta}{Rr_k}]$. Set
\begin{equation}\label{define annulus}
\ {\mathcal A^i_k}:=B_{2^{i+1}Rr_k}\setminus B_{2^iR r_k}\ 
\ {\rm{and}}\ \ {\mathcal B^i_k}:=B_{2^{i+2}Rr_k}\setminus B_{2^{i-1} Rr_k}, 
\ 1\le i\le N_k-1.
\end{equation}
Now we define a radial biharmonic function $v_k$ on the annulus $B_{2\delta}\backslash B_{Rr_k}$ as follows. For simplicity, we
may assume $\frac{2\delta}{Rr_k}$ is a positive integer so that 
$\displaystyle B_{2\delta}\backslash B_{Rr_k}=\bigcup_{i=0}^{N_k-1}
\mathcal A^i_k$. For $0\le i\le N_k-1$, 
$v_k(x)=v_k(|x|)$ satisfies 
\begin{equation}
\begin{cases}\Delta^2 v_k=0 & \mbox{in}\ \ \mathcal A^i_k,\\
 v_k(r)= -\!\!\!\!\!\!\int_{\partial B_{2^{i+1}Rr_k}}u_k,\quad   v'_k(r)= -\!\!\!\!\!\!\int_{\partial B_{2^{i+1}Rr_k}}\frac {\partial u_k}{\partial r},\ &\mbox{if}\  r=2^{i+1}Rr_k,\\
 v_k(r)= -\!\!\!\!\!\!\int_{\partial B_{2^{i}Rr_k}}u_k,\qquad  v'_k(r)= -\!\!\!\!\!\!\int_{\partial B_{2^{i}Rr_k}}\frac {\partial u_k}{\partial r},\  & \mbox{if}\  r=2^{i}Rr_k.
\end{cases}
\end{equation}
Here $-\!\!\!\!\!\!\int$ denotes the average integral. By the standard estimate of biharmonic functions (see, e.g., \cite{HM} Lemma 5.1)
and (\ref{holder annulus}), we have that
$v_k\in W^{4,p}(\mathcal A_k^i)$ for $0\le i\le N_k-1$ and
$$
\Big [v_k\Big]_{C^{\alpha}(\mathcal A^i_k)}\le C\Big(\Big[u_k\Big]_{C^{\alpha}(\mathcal A^i_k)}
+\Big[\nabla u_k\Big]_{L^\infty(\mathcal A_k^i)}\Big)\le C\varepsilon.
$$
In particular, we have
\begin{equation}\label{oscillation-annulus}
\sup_{0\le i\le N_k-1}\mbox{osc}_{\mathcal A^i_k}(u_k-v_k)\le C\varepsilon.
\end{equation}

We now perform the estimate, similar to the arguments by \cite{SaU} or \cite{DT} on harmonic maps and \cite{HM} on biharmonic maps. First, since $u_k-v_k\in W^{4,p}(\mathcal A_k^i)$, we can apply the Green's identity to get that for $0\le i\le N_k-1$,
\begin{eqnarray}\label{Green identity}
\int_{\mathcal A^i_k} \Delta^2 (u_k-v_k)(u_k-v_k)&=&\int_{\mathcal A^i_k} |\Delta (u_k-v_k)|^2+\int_{\partial \mathcal A^i_k} \frac {\partial(\Delta(u_k-v_k))}{\partial \nu}(u_k-v_k)\nonumber\\
&-&\int_{\partial \mathcal A^i_k} \frac {\partial(u_k-v_k)}{\partial \nu}\Delta(u_k-v_k).
\end{eqnarray}
Summing over $0\le i\le N_k-1$, we obtain
\begin{eqnarray}\label{sum over}
\int_{B_{2\delta}\backslash B_{Rr_k}} |\Delta (u_k-v_k)|^2&=&\sum_{i=0}^{N_k-1}\int_{\mathcal A^i_k} \Delta^2 (u_k-v_k)(u_k-v_k)\nonumber\\
&&+\Big(\int_{\partial B_{2\delta}}-\int_{\partial B_{Rr_k}}\Big )\frac {\partial(u_k-v_k)}{\partial \nu}\Delta(u_k-v_k)\nonumber\\
&&-\Big(\int_{\partial B_{2\delta}}-\int_{\partial B_{Rr_k}}\Big )\frac {\partial(\Delta(u_k-v_k))}{\partial \nu}(u_k-v_k)\nonumber\\
&=&\sum_{i=0}^{N_k-1}\int_{\mathcal A^i_k} \Delta^2 u_k (u_k-v_k)
+\Big(\int_{\partial B_{2\delta}}-\int_{\partial B_{Rr_k}}\Big )\frac {\partial(u_k-v_k)}{\partial \nu}\Delta u_k\nonumber\\
&&-\Big(\int_{\partial B_{2\delta}}-\int_{\partial B_{Rr_k}}\Big )\frac {\partial\Delta u_k}{\partial \nu}(u_k-v_k).
\end{eqnarray}
Here we haved use that $\Delta^2 v_k=0$ in $\mathcal A_k^i$, and 
$\displaystyle\int_{\partial B_\rho} \frac{\partial (u_k-v_k)}{\partial\nu}\Delta v_k
=\int_{\partial B_\rho}\frac{\partial \Delta v_k}{\partial\nu} (u_k-v_k)=0$
for $\rho=2\delta$ and $Rr_k$ due to the radial form of $v_k$ and the choices of
boundary conditions of $v_k$.

We can check that the last two terms in the right hand side of (\ref{sum over})
converge to zero as $k\rightarrow \infty$. In fact, by (\ref{Fubini0}), H\"older inequality, and
(\ref{small_energy})
we have that
\begin{eqnarray}\label{1-1-boundary-est}
\Big|\int_{\partial B_{2\delta}}\frac {\partial(u_k-v_k)}{\partial \nu}\Delta u_k \Big|
&\le& \int_{\partial B_{2\delta}} |\nabla u_k||\Delta u_k|
+ (-\!\!\!\!\!\!\int_{\partial B_{2\delta}}|\nabla u_k|)\int_{\partial B_{2\delta}} |\Delta u_k|\nonumber\\
&\le & C\Big(\delta\int_{\partial B_{2\delta}}|\nabla u_k|^4\Big)^\frac14
\Big(\delta\int_{\partial B_{2\delta}}|\nabla u_k|^2\Big)^\frac12\le C\varepsilon^{\frac 32}.
\end{eqnarray}
Similarly, 
\begin{eqnarray}\label{R-boundary-est}
\Big|\int_{\partial B_{Rr_k}}\frac {\partial(u_k-v_k)}{\partial \nu}\Delta u_k\Big|
&\le&   \int_{\partial B_{Rr_k}} |\nabla u_k||\Delta u_k|
+ (-\!\!\!\!\!\!\int_{\partial B_{Rr_k}}|\nabla u_k|)\int_{\partial B_{Rr_k}} |\Delta u_k|\nonumber\\
&\le& C\Big(Rr_k\int_{\partial B_{Rr_k}}|\nabla u_k|^4\Big)^\frac14
\Big(Rr_k\int_{\partial B_{Rr_k}}|\nabla u_k|^2\Big)^\frac12\le C\varepsilon^{\frac 32}.
\end{eqnarray}
For the last term, by Fubini's theorem and (\ref{oscillation-annulus}) we have
\begin{eqnarray}\label{1-2-boundary-est}
\Big|\int_{\partial B_{2\delta}}\frac {\partial\Delta u_k}{\partial \nu}(u_k-v_k)\Big|
&\le &C\Big(\delta\int_{\partial B_{2\delta}}|\nabla^3 u_k|^\frac43\Big)^\frac34\max_{\partial B_{2\delta}}|u_k-v_k|\nonumber\\
&\le& C\epsilon \Big(\int_{B_{4\delta}\setminus B_{\delta}}
|\nabla^3u_k|^\frac43\Big)^\frac34
\le C\varepsilon.
\end{eqnarray}
and, similarly,
\begin{eqnarray}
\Big|\int_{\partial B_{Rr_k}}\frac {\partial\Delta u_k}{\partial \nu}(u_k-v_k)\Big|
&\le& C\Big(Rr_k\int_{\partial B_{Rr_k}}|\nabla^3 u_k|^\frac43\Big)^\frac34\max_{\partial B_{Rr_k}}|u_k-v_k|\nonumber\\
&\le& C\varepsilon \Big(\int_{B_{2Rr_k}\setminus B_{\frac12 Rr_k}}|\nabla^3 u_k|^\frac43\Big)^\frac34\le C\varepsilon.
\end{eqnarray}
Therefore we conclude that the last two terms in the right hand side of (\ref{sum over}) converge to zero as $k\rightarrow \infty$.

For the first term in the right hand side of (\ref{sum over}), we proceed as follows. First, we
can rewrite the equation for $u_k$ as
\begin{equation}\label{fourth-eqn}
\Delta^2 u_k =\mbox{div}(E(u_k))+G(u_k)+h_k,
\end{equation}
where
\begin{equation}\label{eqn-detail}
|E(u_k)|\le C\Big(|\nabla^2u_k||\nabla u_k|+|\nabla u_k|^3 \Big), \quad
|G(u_k)|\le C\Big(|\nabla^2 u_k|^2+|\nabla u_k|^4 \Big).
\end{equation}
Hence
\begin{eqnarray}\label{first-term}
\int_{\mathcal A^i_k} \Delta^2 u_k (u_k-v_k)&=&
\int_{\mathcal A^i_k} \mbox{div} (E(u_k)) (u_k-v_k)
+\int_{\mathcal A^i_k}G(u_k)(u_k-v_k)+\int_{\mathcal A^i_k}h_k(u_k-v_k)\nonumber\\
&=& I^i_k+II^i_k+III^i_k.
\end{eqnarray}
By (\ref{oscillation-annulus}) we have
\begin{equation}\label{II-esti}
|II^i_k|\le C \Big(\int_{\mathcal A^i_k} |\nabla^2 u_k|^2+|\nabla u_k|^4 \Big)\mbox{osc}_{\mathcal A^i_k}u_k \le C\varepsilon \int_{\mathcal A^i_k} |\nabla^2 u_k|^2+|\nabla u_k|^4
\end{equation}
and
\begin{equation}\label{III-esti}
|III^i_k|\le C \mbox{osc}_{\mathcal A^i_k}u_k \int_{\mathcal A^i_k}|h_k| \le C\varepsilon \int_{\mathcal A^i_k}|h_k|.
\end{equation}
For $I^i_k$, by integration by parts we obtain
\begin{eqnarray}\label{I-esti-1}
I^i_k=\int_{\mathcal A^i_k} E(u_k)\cdot \nabla (u_k-v_k)
+\int_{\partial\mathcal A^i_k} E(u_k)(u_k-v_k)\cdot \nu,
\end{eqnarray}
so that after summing over $0\le i\le N_k-1$ we have
\begin{eqnarray}\label{I-esti-2}
\sum_{i=0}^{N_k-1}I^i_k&=&\sum_{i=0}^{N_k-1}\int_{\mathcal A^i_k} E(u_k)\cdot\nabla (u_k-v_k)+
\Big(\int_{\partial B_{2\delta}}-\int_{\partial B_{Rr_k}}\Big)E(u_k)(u_k-v_k)\cdot \nu\nonumber\\
&\le & C\sum_{i=0}^{N_k-1}\int_{\mathcal A^i_k}\Big(|\nabla^2u_k||\nabla u_k|+|\nabla u_k|^3 \Big)|\nabla (u_k-v_k)|\nonumber\\
&+&C\Big(\int_{\partial B_{2\delta}}+\int_{\partial B_{Rr_k}}\Big) \Big(|\nabla^2u_k||\nabla u_k|+|\nabla u_k|^3 \Big)|u_k-v_k|\nonumber\\
&\le & IV_k+V_k.
\end{eqnarray}
By (\ref{Fubini0}) and H\"older inequality, we see that
$$V_k\le C\varepsilon.$$
For $IV_k$, we proceed the estimate as follows.
By Nirenberg's interpolation inequality and (\ref{oscillation-annulus}), we have
\begin{eqnarray}\label {Nirenberg}
\Big(\int_{\mathcal A^i_k}|\nabla u_k|^4 \Big)^{\frac 14}
&\le& C\Big[u_k\Big]^{\frac 12}_{{\rm{BMO}}(\mathcal B^i_k)}\Big(\int_{\mathcal B^i_k}(|\nabla^2 u_k|^2+\frac{|\nabla u_k|^2}{(2^iRr_k)^2})\Big)^{\frac 14}\nonumber\\
&\le& C\varepsilon^\frac12\Big(\int_{\mathcal B^i_k}(|\nabla^2 u_k|^2+\frac{|\nabla u_k|^2}{(2^iRr_k)^2})\Big)^{\frac 14}.
\end{eqnarray}
Since $v_k$ is a biharmonic function on $\mathcal A_k^i$, we also have
\begin{eqnarray}\label {Nirenberg-ineq}
\Big(\int_{\mathcal A^i_k}|\nabla(u_k-v_k)|^4 \Big)^{\frac 14}
&\le& C\Big(\int_{\mathcal A^i_k}|\nabla u_k |^4 \Big)^{\frac 14}\nonumber\\
&\le& C\Big[u_k\Big]^{\frac 12}_{{\rm{BMO}}(\mathcal B^i_k)}\Big(\int_{\mathcal B^i_k}|\nabla^2 u_k|^2+\frac{|\nabla u_k|^2}{(2^iRr_k)^2}\Big)^{\frac 14}\nonumber\\
&\le& C\varepsilon^\frac12\Big(\int_{\mathcal B^i_k}(|\nabla^2 u_k|^2+\frac{|\nabla u_k|^2}{(2^iRr_k)^2})\Big)^{\frac 14}.
\end{eqnarray}
Therefore we have
\begin{eqnarray}\label{IV-1-esti}
IV_k&=&\sum_{i=0}^{N_k-1}\int_{\mathcal A^i_k}|\nabla^2u_k||\nabla u_k||\nabla (u_k-v_k)|\nonumber\\
&\le &
C\sum_{i=0}^{N_k-1}\Big(\int_{\mathcal A^i_k}|\nabla^2 u_k|^2\Big)^{\frac 12}\Big(\int_{\mathcal A^i_k}|\nabla u_k|^4 \Big)^{\frac 14}
\Big(\int_{\mathcal A^i_k}|\nabla (u_k-v_k)|^4\Big)^{\frac 14}\nonumber\\
&\le & C\varepsilon \Big(\sum_{i=0}^{N_k-1}\int_{\mathcal A^i_k}|\nabla^2 u_k|^2\Big)^{\frac 12}\Big(\sum_{i=0}^{N_k-1}\int_{\mathcal B^i_k}|\nabla^2 u_k|^2+\frac{|\nabla u_k|^2}{(2^iRr_k)^2}\Big)^{\frac 12}\nonumber\\
&\le & C\varepsilon
\Big(\int_{B_{2\delta}\backslash B_{Rr_k} }|\nabla^2 u_k|^2\Big)^{\frac 12}
\Big(\int_{B_{4\delta}\backslash B_{\frac12Rr_k} }|\nabla^2 u_k|^2+\frac{|\nabla u_k|^2}
{|x|^2}\Big)^{\frac 12}.
\end{eqnarray}
Applying Lemma 5.2 in \cite{HM}, we have the following Hardy inequality:
\begin{equation}\label{Hardy-ineq}
\int_{B_{4\delta}\backslash B_{\frac12Rr_k}} \frac{|\nabla u_k|^2}{|x|^2}\le
\int_{B_{4\delta}\backslash B_{\frac12 Rr_k}} |\nabla^2 u_k|^2+\Big(\int_{\partial B_{4\delta}}-\int_{\partial B_{\frac12Rr_k}}\Big)\big(
\frac 1{|x|}|\nabla u_k|^2\big).
\end{equation}
By (\ref{Fubini0}) and H\"older inequality, we have that
\begin{eqnarray}\label{Fubini1}
\frac 1{\delta}\int_{\partial B_{4\delta}}|\nabla u_k|^2 \le  C
\Big(\delta\int_{\partial B_{4\delta}}|\nabla u_k|^4\Big)^\frac12
\le C\Big(\int_{B_{8\delta}\backslash B_{2\delta}}|\nabla u_k|^4\Big)^{\frac 12}\le C\varepsilon,
\end{eqnarray}
and, similarly,
\begin{eqnarray}\label{Fubini2}
\frac 1{Rr_k}\int_{\partial B_{\frac12 Rr_k}}|\nabla u_k|^2
\le C\Big(Rr_k\int_{\partial B_{\frac12 Rr_k}}|\nabla u_k|^4\Big)^\frac12
\le C \Big(\int_{B_{Rr_k}\setminus B_{\frac14 Rr_k}}|\nabla u_k|^4\Big)^\frac12
\le C\varepsilon.
\end{eqnarray}
Substituting (\ref{Hardy-ineq}), (\ref{Fubini1}) and (\ref{Fubini2}) into (\ref{IV-1-esti}), 
we obtain 
\begin{eqnarray}\label{IV-1-esti-2}
IV_k\le C\varepsilon\Big( \int_{B_{4\delta}\backslash B_{\frac12Rr_k}} |\nabla^2 u_k|^2\Big)
+C\varepsilon.
\end{eqnarray}
Using the same argument to estimate the second term of $IV_k$, we have
\begin{eqnarray}\label{IV-summary}
IV_k\le  C\varepsilon\Big( \int_{B_{2\delta}\backslash B_{Rr_k}} |\nabla^2 u_k|^2\Big)
+C\varepsilon.
\end{eqnarray}
Combining the estimates together yields 
\begin{equation}\label{angular-energy}
\int_{B_{2\delta}\backslash B_{Rr_k}} |\Delta (u_k-v_k)|^2\le C\varepsilon\Big(\int_{B_{4\delta}\backslash B_{\frac12Rr_k}}
|\nabla^2 u_k|^2+|\nabla u_k |^4+|h_k|\Big)+C\varepsilon.
\end{equation}
This, combined with Calderon-Zygmund's $W^{2,2}$ estimate,
yields
\begin{equation}\label{angular-energy1}
\int_{B_{\delta}\backslash B_{2Rr_k}} |\nabla^2 (u_k-v_k)|^2\le C\varepsilon\Big(\int_{B_{4\delta}\backslash B_{\frac12Rr_k}}
|\nabla^2 u_k|^2+|\nabla u_k |^4+|h_k|\Big)+C\varepsilon.
\end{equation}

\noindent{\bf Step 2}. {\it Control of radial component of hessian energy in the neck region}:\\

Since $v_k$ is radial, it is easy to see that (\ref{angular-energy1})
yields
\begin{eqnarray}\label{angular-energy2}
\int_{B_\delta\setminus B_{2Rr_k}}|\nabla_T\nabla u_k|^2
&\le& C\varepsilon\Big(\int_{B_{4\delta}\backslash B_{\frac12Rr_k}}
|\nabla^2 u_k|^2+|\nabla u_k |^4+|h_k|\Big)+C\varepsilon\nonumber\\
&\le& C\varepsilon.
\end{eqnarray}
Here $\displaystyle\nabla_T\nabla u_k=\nabla^2u_k-\frac{\partial}{\partial r}(\nabla u_k)$ denotes the tangential component of $\nabla u_k$. 

Next, we want to apply the Pohozaev type argument to $W^{4,p}$-approximate biharmonic maps $u_k$ with $L^p$ bi-tension field $h_k$ for $p\ge\frac43$ to control
$\displaystyle\int_{B_\delta\setminus B_{2Rr_k}}\Big|\frac{\partial^2u_k}{\partial r^2}\Big|^2$
by $\displaystyle\int_{B_\delta\setminus B_{2Rr_k}}|\nabla_T \nabla u_k|^2$
and $\|h_k\|_{L^p(B_{2\delta})}$.
This type of argument is well-known in the blow up analysis of harmonic or approximate
harmonic maps on Riemann surfaces (see \cite{SaU}, \cite{LW1}, \cite{LW}, and 
\cite{LinR}). In the context of biharmonic maps, it was first
derived by Hornung-Moser \cite{HM}. 

By (\ref{small_energy}) and Lemma 2.3, we see that $u_k\in W^{4,p}(B_{2\delta}\setminus B_{\frac12 Rr_k})$. While in $B_{\frac12Rr_k}$,
since $\widetilde{u}_k(x)=u_k(r_kx): B_R\to N$ is an approximate biharmonic map that converges to the bubble $\omega_1$,
we conclude that $\displaystyle\|\widetilde{u}_k\|_{W^{4,p}(B_R)}<+\infty$ and hence $u_k\in W^{4,p}(B_{\frac12Rr_k})$
by scaling. From this,  we then see 
$u_k\in W^{4,p}(B_{2\delta})$. Since $x\cdot\nabla u_k\in L^4(B_{2\delta})$ and $p\ge\frac43$, we see that 
$\Delta^2 u_k \cdot (x\cdot\nabla u_k)\in L^1(B_{2\delta})$
and $h_k\cdot (x\cdot\nabla u_k)\in L^1(B_{2\delta})$. 
Since the equation (\ref{approx_biharm}) 
implies that $(\Delta^2u_k-h_k)(x)\perp T_{u_k(x)}N$ a.e. $x\in B_{2\delta}$. Note also  that
$x\cdot\nabla u_k(x)\in T_{u_k(x)}N$ for a.e. $x\in B_{2\delta}$. 
Multiplying the equation (\ref{approx_biharm}) by $x\cdot \nabla u_k$ and integrating over $B_r$
for any $0<r\le 2\delta$, we have
\begin{eqnarray}\label{id1}
\int_{B_r}\Delta^2u_k\cdot(x\cdot \nabla u_k)=\int_{B_r}h_k\cdot(x\cdot \nabla u_k).
\end{eqnarray}
Applying Green's identity, we have
\begin{eqnarray}\label{LHS-identity}
&&\int_{B_r}\Delta^2u_k\cdot(x\cdot \nabla u_k)\nonumber\\
&=&
\int_{B_r}\Delta u_k \Delta(x\cdot\nabla u_k)
+r\int_{\partial B_r}\frac{\partial}{\partial r}(\Delta u_k)\frac{\partial u_k}{\partial r}
-\int_{\partial B_r}\Delta u_k(\frac{\partial u_k}{\partial r}+r\frac{\partial^2 u_k}{\partial r^2}).
\end{eqnarray}
Direct calculations yield
\begin{eqnarray}\label{LHS-identity1}
\int_{B_r}\Delta u_k \Delta(x\cdot\nabla u_k)
&=& \int_{B_r} x\cdot\nabla\Big(\frac{|\Delta u_k|^2}2\Big)+2|\Delta u_k|^2\nonumber\\
&=& \int_{B_r} {\rm{div}}\Big(\frac{|\Delta u_k|^2}2 x\Big)
=r\int_{\partial B_r} \frac{|\Delta u_k|^2}2.
\end{eqnarray}
Putting (\ref{LHS-identity1}),  (\ref{LHS-identity}), and (\ref{id1})  together yields
\begin{equation}\label{id2}
r\int_{\partial B_r} \frac{|\Delta u_k|^2}2
+r\int_{\partial B_r}\frac{\partial}{\partial r}(\Delta u_k)\frac{\partial u_k}{\partial r}
-\int_{\partial B_r}\Delta u_k(\frac{\partial u_k}{\partial r}+r\frac{\partial^2 u_k}{\partial r^2})
=\int_{B_r}h_k\cdot(x\cdot \nabla u_k).
\end{equation}
Applying integration by parts multi-times to (\ref{id2}) in the same way as 
\cite{CWY} or Angelsberg \cite{an}, we can obtain that for a.e. $0<r\le 2\delta$,
\begin{eqnarray}\label{id3}
\int_{\partial B_r}|\Delta u_k|^2
&=&4\int_{\partial B_r} \Big(\frac{|u_{k,\alpha}
+x^\beta u_{k,\alpha\beta}|^2}{r^2}
+2\frac{|x\cdot\nabla u_k|^2}{r^4}\Big)\nonumber\\
&+&2\frac{d}{dr}\int_{\partial B_r}
\Big(-\frac{x^\alpha u_{k,\beta} u_{k,\alpha\beta}}{r}
+2\frac{|x\cdot\nabla u_k|^2}{r^3}-2\frac{|\nabla u_k|^2}{r}\Big)\nonumber\\
&+&\frac{1}{r}\int_{B_r} h_k\cdot(x\cdot\nabla u_k).
\end{eqnarray}
Recall that in the spherical coordinate, we have
$$\Delta u_k=u_{k,rr}+\frac{3}{r} u_{k,r}+\frac{1}{r^2}\Delta_{S^3}u_k,$$
where $\Delta_{S^3}$ denotes the Laplace operator on the standard three sphere $S^3$.
Hence we have
\begin{eqnarray}\int_{\partial B_r}|\Delta u_k|^2
&=&\int_{\partial B_r}
\Big[|u_{k,rr}|^2+\frac{9}{r^2}|u_{k,r}|^2+\frac{6}{r} u_{k,r}u_{k,rr}\Big]\nonumber\\
&+&\int_{\partial B_r}\Big[\frac{1}{r^4}|\Delta_{S^3} u_k|^2+(u_{k,rr}+\frac{3}{r}u_{k,r})\cdot
(\frac{2}{r^2}\Delta_{S^3}u_k)\Big].\label{id4}
\end{eqnarray}
On the other hand, we have
\begin{eqnarray}
&&4\int_{\partial B_r} \Big(\frac{|u_{k,\alpha}
+x^\beta u_{k,\alpha\beta}|^2}{r^2}
+2\frac{|x\cdot\nabla u_k|^2}{r^4}\Big)\nonumber\\
&=&\int_{\partial B_r}\Big[\frac{12}{r^2}|u_{k,r}|^2
+4|u_{k,rr}|^2+\frac{4}{r^2}|\nabla_{S^3} u_{k,r}|^2
+\frac{8}{r}u_{k,r} u_{k,rr}\Big].\label{id5}
\end{eqnarray}
Substituting (\ref{id4}) and (\ref{id5}) into (\ref{id3}) and integrating over
$r\in [Rr_k, \delta]$, we obtain
\begin{eqnarray}
&&\int_{B_\delta\setminus B_{Rr_k}}
\Big[3|u_{k,rr}|^2+\frac{3}{r^2}|u_{k,r}|^2+\frac{2}{r} u_{k,r} u_{k,rr}\Big]\nonumber\\
&\le& \int_{B_\delta\setminus B_{Rr_k}} 
\Big[ \frac{1}{r^4}|\Delta_{S^3} u_k|^2+(u_{k,rr}+\frac{3}{r}u_{k,r})\cdot
(\frac{2}{r^2}\Delta_{S^3}u_k)\Big]+\int_{Rr_k}^{\delta} \int_{B_r}|h_k||u_{k,r}|\nonumber\\
&+& 2\Big(\int_{\partial B_\delta}-\int_{\partial B_{Rr_k}}\Big)
\Big(\frac{x^\alpha u_{k,\beta} u_{k,\alpha\beta}}{r}
-2\frac{|x\cdot\nabla u_k|^2}{r^3}+2\frac{|\nabla u_k|^2}{r}\Big)\nonumber\\
&=& I_k+II_k+III_k.
\end{eqnarray}
By H\"older inequality, (\ref{angular-energy2}),  (\ref{Hardy-ineq}),
(\ref{Fubini1}), and (\ref{Fubini2}),  we have
$$|I_k|\le \int_{B_\delta\setminus B_{Rr_k}}|\nabla_T\nabla u_k|^2
+\Big(\int_{B_\delta\setminus B_{Rr_k}}|\nabla^2 u_k|^2\Big)^\frac12
\Big(\int_{B_\delta\setminus B_{Rr_k}}|\nabla_T\nabla u_k|^2\Big)^\frac12+C\varepsilon
\le C\varepsilon.$$
For $II_k$ we have
$$|II_k|
\le \delta \|h_k\|_{L^p(B_\delta)}\|\nabla u_k\|_{L^{\frac{p}{p-1}}(B_\delta)}
\le C\delta \|\nabla u_k\|_{L^4(B_\delta)}\le C\delta,$$
where we have used the fact $p\ge\frac43$ so that $\frac{p}{p-1}\le 4$.

We use (\ref{Fubini0}) to estimate $III_k$ as follows. First we have
\begin{eqnarray*}
&&\Big|\int_{\partial B_\delta}\Big(\frac{x^\alpha u_{k,\beta} u_{k,\alpha\beta}}{r}
-2\frac{|x\cdot\nabla u_k|^2}{r^3}+2\frac{|\nabla u_k|^2}{r}\Big)\Big|\\
&\le& C\Big[\int_{\partial B_\delta} |\nabla u_k||\nabla^2u_k|
+\delta^{-1}\int_{\partial B_\delta}|\nabla u_k|^2\Big]\\
&\le& C\Big(\delta\int_{\partial B_\delta}|\nabla u_k|^4\Big)^\frac14
\Big(\delta\int_{\partial B_\delta}|\nabla^2 u_k|^2\Big)^\frac12+C\Big(\delta\int_{\partial B_\delta}|\nabla u_k|^4\Big)^\frac12\\
&\le& C\Big[\|\nabla u_k\|_{L^4(B_{2\delta}\setminus B_{\frac{\delta}2})}
\|\nabla^2 u_k\|_{L^2(B_{2\delta}\setminus B_{\frac{\delta}2})}
+\|\nabla u_k\|_{L^4(B_{2\delta}\setminus B_{\frac{\delta}2})}^2\Big]
\le C\varepsilon.
\end{eqnarray*}
Similarly, by (\ref{Fubini0}) we have
\begin{eqnarray*}
&&\Big|\int_{\partial B_{Rr_k}}\Big(\frac{x^\alpha u_{k,\beta} u_{k,\alpha\beta}}{r}
-2\frac{|x\cdot\nabla u_k|^2}{r^3}+2\frac{|\nabla u_k|^2}{r}\Big)\Big|\\
&\le& C\Big[\int_{\partial B_{Rr_k}} |\nabla u_k||\nabla^2u_k|
+(Rr_k)^{-1}\int_{\partial B_{Rr_k}}|\nabla u_k|^2\Big]\\
&\le& 
C\Big(Rr_k\int_{\partial B_{Rr_k}} |\nabla u_k|^4\Big)^\frac14 \Big(Rr_k\int_{\partial B_{Rr_k}} |\nabla^2 u_k|^2\Big)^\frac12
+C\Big(Rr_k\int_{\partial B_{Rr_k}} |\nabla u_k|^4\Big)^\frac12\\
&\le& C\Big[\|\nabla u_k\|_{L^4(B_{2Rr_k}\setminus B_{\frac12Rr_k})}
\|\nabla^2 u_k\|_{L^2(B_{2Rr_k}\setminus B_{\frac12Rr_k})}
+\|\nabla u_k\|_{L^4(B_{2Rr_k}\setminus B_{\frac12Rr_k})}^2\Big]
\le C\varepsilon.
\end{eqnarray*}
Therefore, by putting these estimates together, we have
\begin{equation}\label{radial-energy}
\int_{B_\delta\setminus B_{Rr_k}}
\Big[3|u_{k,rr}|^2+\frac{3}{r^2}|u_{k,r}|^2+\frac{2}{r} u_{k,r} u_{k,rr}\Big]
\le C(\varepsilon+\delta).
\end{equation}
\bigskip
\noindent{\bf Competition of Proof of Theorem \ref{energy_identity}}: \\

Since 
$\displaystyle \frac{2}{r} u_{k,r} u_{k,rr}\ge -(|u_{k,rr}|^2+\frac{1}{r^2}|u_{k,r}|^2)$,
(\ref{radial-energy}) implies
$$\int_{B_\delta\setminus B_{Rr_k}} |u_{k,rr}|^2
\le C(\varepsilon+\delta),$$
this, combined with (\ref{angular-energy2}), implies
$$\int_{B_\delta\setminus B_{Rr_k}}|\nabla^2 u_k|^2
\le C(\varepsilon+\delta).$$
Thus there is no  concentration of hessian energy in the neck region. It is
well known that this yields the energy identity (\ref{energy_id1}). 
To show (\ref{energy_id2}),  observe that
Nirenberg's interpolation inequality and (\ref{energy_id1}) imply
\begin{eqnarray*}
\|\nabla u_k\|_{L^4(B_\delta\setminus B_{Rr_k})}^2
&\le& C\|\nabla u_k\|_{L^\infty(B_{2\delta})}
\Big(\|\nabla u_k\|_{L^2(B_{2\delta})}
+\|\nabla^2 u_k\|_{L^2(B_{2\delta}\setminus B_{\frac12 Rr_k})}\Big)\\
&\le& C(\epsilon+\delta+o(1))
\end{eqnarray*}
where we have used that
$\displaystyle \|\nabla u_k\|_{L^2(B_{2\delta})}=\|\nabla u\|_{L^2(B_{2\delta})}+o(1)
=o(1).$ Thus (\ref{energy_id2}) also holds.
\qed\\

\noindent{\bf Proof of Corollary \ref{energy_identity1}}:

\smallskip
It follows from the energy inequality (1.8)
that there exists $t_k\uparrow +\infty$ such that
$u_k(\cdot)=u(\cdot,t_k)$ is an approximate biharmonic map into
$N$ with bi-tension field
$h_k=u_t(\cdot, t_k)\in L^2(\Omega)$ satisfying
$$\Big\|h_k\Big\|_{L^2(\Omega)}
=\Big\|u_t(\cdot,t_k)\Big\|_{L^2(\Omega)}\rightarrow 0.$$
Moreover,
$$\Big\|u_k\Big\|_{W^{2,2}(\Omega)}\le
C\Big\| u_0\Big\|_{W^{2,2}(\Omega)}.$$
Therefore we may assume that after taking another subsequence,
$u_k\rightharpoonup u_\infty$ in $W^{2,2}(\Omega,N)$.
It is easy to see that $u_\infty$ is a biharmonic map so that
$u_\infty\in C^\infty(\Omega,N)$ (see \cite{W2}).
All other conclusions follow directly from Theorem \ref{energy_identity}. \qed


\begin{thebibliography}{10}

\bibitem {A} Adams, D. R.: A note on Riesz potentials. {\em Duke Math. J.}  {\bf 42} (4), 765-778 (1975).

\bibitem{an} G. Angelsberg, A monotonicity formula for stationary biharmonic maps. 
{\em Math. Z.},  {\bf 252} (2006), 287-293.




\bibitem {CWY} Chang, A., Wang, L., Yang, P.: A
regularity theory of biharmonic maps. {\em Comm. Pure Appl. Math.} {\bf 52},
1113-1137 (1999).


\bibitem {DT} Ding, W. Y., Tian, G.: Energy identity for a class of approximate
harmonic maps from surfaces. {\em Comm. Anal. Geom.} {\bf 3}, 543-554 (1995).



\bibitem{Gastel} Gastel, A.: The extrinsic polyharmonic map heat flow in the critical dimension.
{\em Adv. Geom.} {\bf  6} no. 4, 501-521 (2006).

\bibitem{GLW} Gong, H. J.,  Lamm, T., Wang, C. Y.:  {Boundary regularity for a class of biharmonic
maps}. {\em  Calc. Var. Partial Differential Equations}, to appear.

\bibitem {H}H\'elein, F.: Harmonic Maps, Conservation Laws, and Moving Frames.
Cambridge Tracts in Mathematics, 150, Cambridge: Cambridge
University Press (2002).


\bibitem{HM} Hornung, P., Moser, R.: Energy identity for instrinsically biharmonic maps in four
dimensions. {\em Anal. PDE}, to appear.

\bibitem {IM}Iwaniec, T., Martin, G.: Quasiregular mappings in even dimensions. {\em Acta Math.}
{\bf 170}, 29-81
(1993).


\bibitem {K} Ku, Y.:  Interior and boundary regularity of
intrinsic biharmonic maps to spheres. {\em Pacific J. Math.} {\bf 234}, 43-67 (2008).

\bibitem{Lamm} Lamm, T.: Heat flow for extrinsic biharmonic maps with small initial energy.
{\em Ann. Global Anal. Geom.}  {\bf 26} no. 4, 369-384 (2004).

\bibitem{LR} Lamm, T., Rivi\`ere, T.: Conservation laws for fourth order systems in four dimensions.
{\em Comm. PDE.} {\bf 33}, 245-262 (2008).



\bibitem {LinR}Lin, F. H., Rivi\`ere, T.: Energy quantization for harmonic maps.
{\em Duke Math. J.} {\bf 111}, 177-193 (2002).

\bibitem {LW1}Lin, F. H., Wang, C. Y.: Energy id entity of harmonic map flows
from surfaces at finite singular time. {\em Calc. Var. Partial
Differential Equations} {\bf 6}, 369-380 (1998).


\bibitem {LW}Lin, F. H., Wang, C. Y.: Harmonic and quasi-harmonic spheres II.
{\em Comm. Anal. Geom.} {\bf10}, 341-375 (2002).

\bibitem{M} Moser, R.: Weak solutions of a biharmonic map heat flow.
{\em Adv. Calc. Var.}  {\bf 2}  no. 1, 73-92 (2009).





\bibitem{scheven1} Scheven, C.:
Dimension reduction for the singular set of biharmonic maps.
{\em Adv. Calc. Var.} {\bf  1} no. 1, 53-91 (2008).


\bibitem{scheven2} Scheven, C.: An optimal partial regularity result for minimizers of an intrinsically defined second-order functional.
{\em Ann. Inst. H. Poincar\'e Anal. Non Lin\'eaire}.  {\bf 26}  no. 5, 1585-1605 (2009).

\bibitem {S}Strzelecki, P.: On biharmonic maps and their generalizations. {\em
 Calc. Var. Partial Differential Equations,} {\bf 18} (4), 401-432 (2003).

\bibitem {SaU} Sacks, J., Uhlenbeck. K.: The existence of minimal immersions
of 2-spheres. Annals of Math., 113, 1-24(1981).




\bibitem{Struwe} Struwe, M.:  Partial regularity for biharmonic maps, revisited.
{\em  Calc. Var. Partial Differential Equations,} {\bf 33} (2), 249-262 (2008).

\bibitem {W1} Wang, C. Y.: Remarks on biharmonic maps into spheres. {\em  Calc. Var. Partial Differential Equations,} {\bf 21}, 221-242 (2004).


\bibitem {W2} Wang, C. Y.: Biharmonic maps from $\mathbb R^4$ into a Riemannian
manifold. {\em Math. Z.} {\bf 247}, 65-87 (2004).

\bibitem{W3} Wang, C. Y.: Stationray biharmonic Maps from $\mathbb R^m$ into a Riemannian Manifold.
{\em Comm. Pure Appl. Math.} {\bf 57}, 0419-0444 (2004).

\bibitem{W4} Wang, C. Y.: Heat flow of biharmonic maps in dimensions four and its application.
{\em Pure Appl. Math. Q.}  {\bf 3 }no. 2, part 1, 595-613 (2007).

\bibitem{WZ} Wang, C. Y., Zheng, S. Z.: Energy identity for a class of approximate biharmonic maps in dimension four. arXiv:1110.6536, preprint (2011).

\end{thebibliography}
\end{document}